\documentclass[a4paper,leqno]{amsart}

\usepackage{amsmath}
\usepackage{amsthm}
\usepackage{amssymb}
\usepackage{amsfonts}
\usepackage[applemac]{inputenc}
\usepackage{mathrsfs}
\usepackage{pdfsync}
\usepackage{graphicx}
\usepackage{mathabx}
\usepackage{color}
\def\R{{\mathbb R}}% real numbers
\def\N{{\mathbb N}}% nonnegative integers
\def\Z{{\mathbb Z}}% integers
\def\le{\leqslant}% lessoreqal
\def\ge{\geqslant}%greaterorequal

\theoremstyle{plain}
\newtheorem{theorem}{Theorem}[section]

\theoremstyle{definition}
\newtheorem{definition}[theorem]{Definition}

\newtheorem{remark}[theorem]{Remark}
\newtheorem*{remark*}{Remark}

\numberwithin{equation}{section}

\begin{document}

\title[Transverse (in-)stability of cubic-quintic solitary waves]
{Numerical study of transverse (in-)stability of solitary waves in 
the cubic-quintic nonlinear Schr\"odinger equation}

\author[C. Klein]{Christian Klein}
\address[C.~Klein]
{Institut de Math\'ematiques de Bourgogne, Universit\'e 
Bourgogne Europe, \\
Institut Universitaire de France,\\
9 avenue Alain Savary, BP 47870, 21078 Dijon Cedex}
\email{christian.klein@u-bourgogne.fr}

\author[C. Sparber]{Christof Sparber}
\address[C.~Sparber]
{Department of Mathematics, Statistics, and Computer Science, M/C 249, University of Illinois at Chicago, 851 S. Morgan Street, Chicago, IL 60607, USA}
\email{sparber@uic.edu}

\begin{abstract}
We study the nonlinear Schr\"odinger equation with a competing cubic-quintic power law nonlinearity on the waveguide domain $\mathbb R_x \times \mathbb T_{L_y}$. 
This model is globally well-posed and admits line solitary wave solutions, whose transverse (in-)stability is numerically investigated. We consider both 
spatially localized perturbations and periodic deformations of the line solitary wave and numerically confirm that there exists a critical torus length $L_y>0$ above which instability appears. 
\end{abstract}

\date{\today}


\thanks{C.K. thanks for support by the ANR project 
ANR-17-EURE-0002 EIPHI and by the ANR project 
ISAAC-ANR-23-CE40-0015-01. C.S. is supported by the MPS Simons foundation through award no. 851720}
\maketitle

%%%%%%%%%%%%%%%%%%%%%%%%%%%%%%%%%%%
%%%%%%%%%%%%%%%%%%%%%%%%%%%%%%%%%%%

\section{Introduction}
\label{sec:intro}

In this work we consider the nonlinear Schr\"odinger equation 
(NLS) with competing cubic-quintic nonlinearity, i.e.,
\begin{equation}\label{nls}
\left\{
\begin{aligned}
  i\partial_t u +\Delta u & =  - |u|^2u+ |u|^4u,\\
   u_{\mid
  t=0} & = u_0.
\end{aligned}  
\right.
\end{equation}
This model appears in numerous applications, including nonlinear optics, superconductivity, 
and the description of Bose-Einstein condensates, see, e.g., \cite{DMM, GFTC, Ma} and the references therein.
In contrast to the case of a single power-law nonlinearity, this model no longer admits a scaling symmetry, which makes 
the mathematical analysis more difficult. Indeed, ever since the seminal paper of Killip et al. \cite{KOPV}, this NLS 
has attracted considerable mathematical interest, cf. \cite{CaSp, FYZ, LeRo, Mu, Oh}. 

Equation \eqref{nls} formally conserves (among others) the mass 
\[
M(u)=\| u(t, \cdot) \|_{L^2}^2
\]
and the total energy
\[
E(u) = \frac{1}{2}\|\nabla u(t, \cdot)\|_{L^2}^2 -\frac{1}{2}\|
u(t, \cdot)\|_{L^4}^4+ \frac{1}{3}\|  u(t, \cdot)\|_{L^6}^6.
\]
Combined with H\"older's inequality 
\[ 
\|u \|_{L^4}^4 \le \|u \|_{L^2}\|u \|_{L^6}^3,
\]
these conservation laws imply that the cubic, focusing part cannot be an obstruction to global
well-posedness for solutions $u(t,\cdot)$ in the energy-space $H^1$. 
This is in sharp contrast to the case without the defocusing quintic nonlinearity, where finite-time blow-up of solutions is well-known in space dimensions two and higher, cf. \cite{CazCourant}. 

\subsection{One- and two-dimensional ground states} A particular class of solutions to \eqref{nls} are solitary waves of the form 
\[
u(t,{\bf x}) = \phi_\omega({\bf x}) e^{i \omega t}, \quad \omega \in \R,
\] 
where $\phi_\omega$ is a localized (i.e., vanishing at infinity) solution to the stationary equation 
 \begin{equation}\label{eq:solitondef}
    -\Delta \phi_\omega + \omega \phi_\omega
    -|\phi_\omega|^2\phi_\omega+|\phi_\omega|^4\phi_\omega=0.
  \end{equation} 
By deriving the associated Pohozaev-Nehari identities, one can prove that 
a necessary condition for the existence of nontrivial profiles $0\not = \phi_\omega \in H^1$ is given by $\omega \in (0, \frac{3}{16})$, see \cite{CaSp}.  

\begin{figure}[htb!]
  \includegraphics[width=0.49\textwidth]{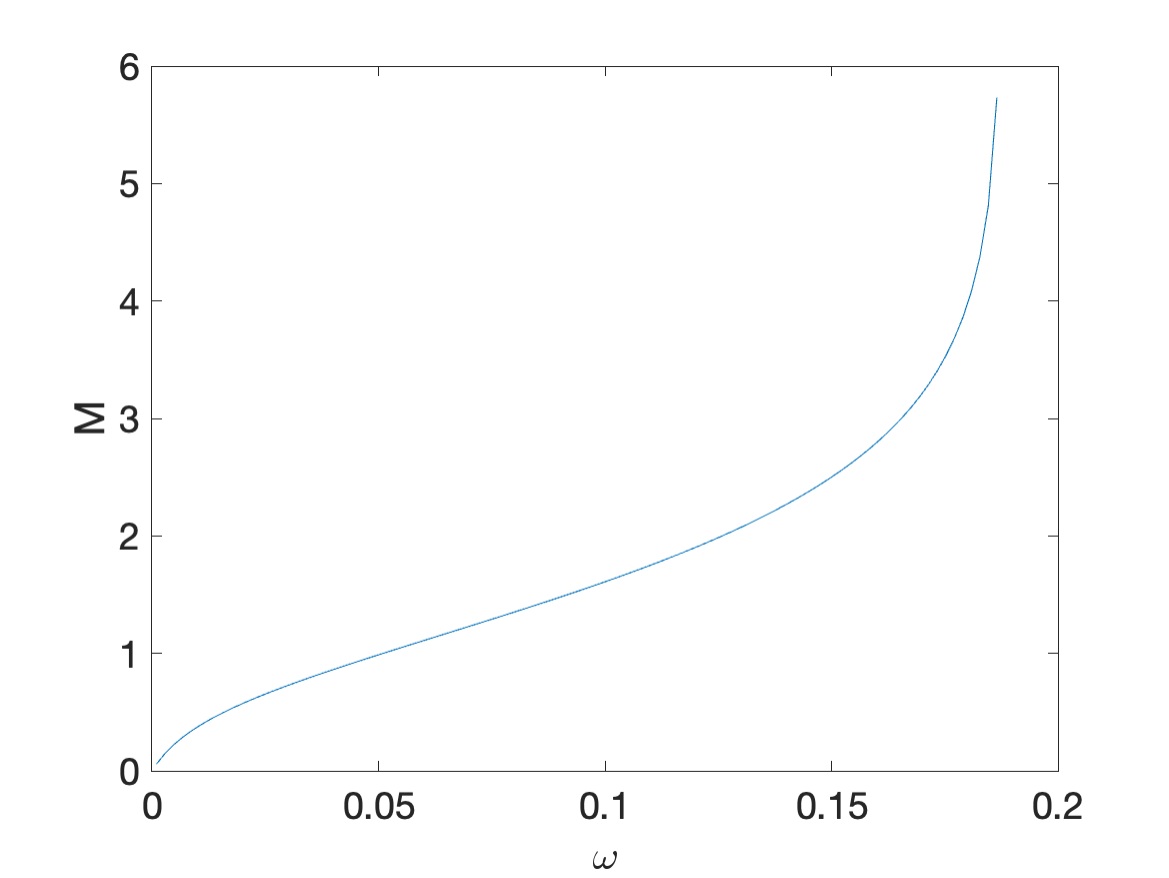}
  \includegraphics[width=0.49\textwidth]{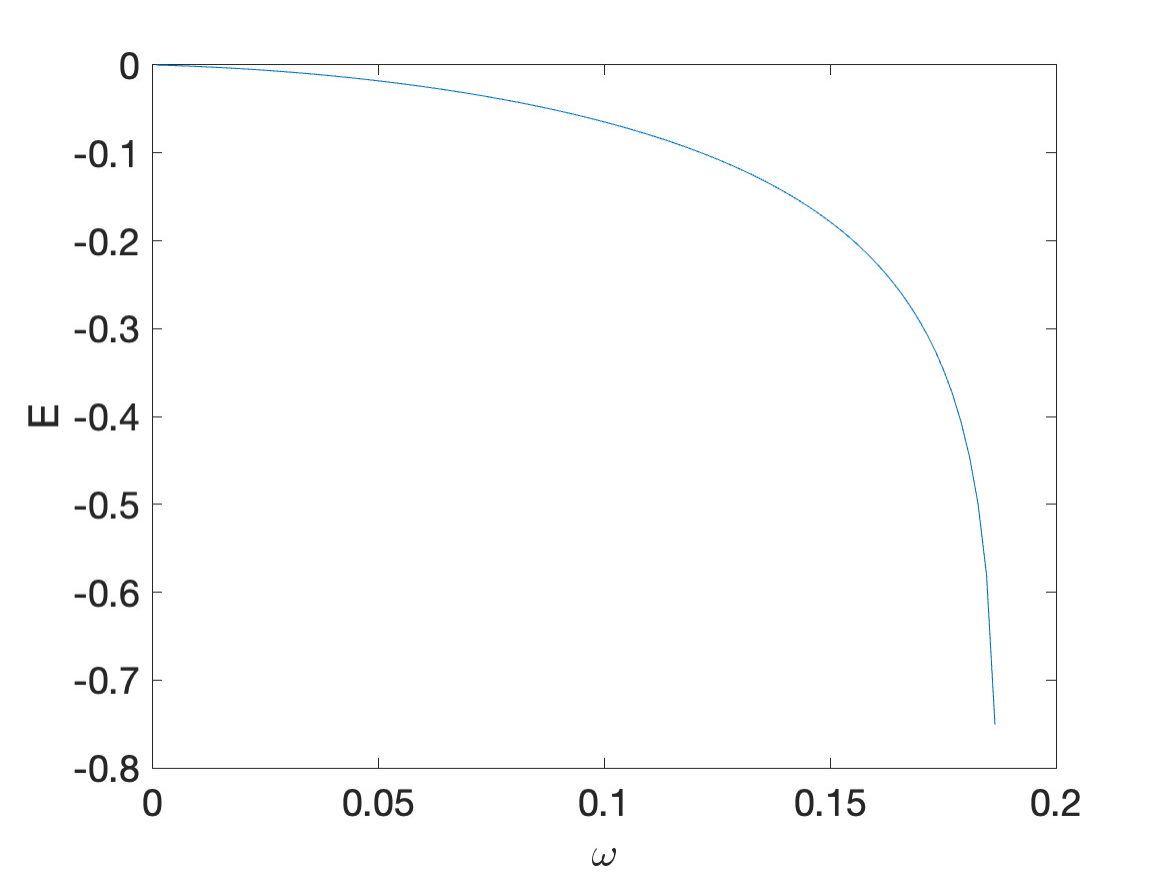}
\caption{The mass (left) and the energy (right) of the 1D solitary 
wave (\ref{phi}) as a function of $\omega$.}
 \label{ME}
\end{figure}

In one space dimension, the profile $\phi_\omega$ can be explicitly computed in the form
\begin{equation}\label{phi}
	    \phi_\omega(x) = 
		\left(\frac{1}{4\omega}+\sqrt{\frac{1}{16\omega^2}-\frac{1}{3\omega}}\cosh\big(2\sqrt{\omega}x\big)\right)^{-1/2}, \quad x \in \R.
\end{equation}
In Fig.~\ref{ME} we show the mass $M(\phi_\omega)$ and the energy $E(\phi_\omega)$ of this explicit family of 1D solitary 
waves as a function of $\omega\in (0, \frac{3}{16})$. It can be seen that both $M(\phi_\omega)\ge 0$ and $E(\phi_\omega)\le 0$ 
vanish as $\omega\to0$, while they both diverge as 
$\omega \to \frac{3}{16}$. Moreover, both are monotonic w.r.t. $\omega$.

Let us briefly recall the definition of orbital stability of solitary waves. The need of the inclusion of the parameters $\theta\in \R$ and ${\bf x}_0\in \R^d$ is thereby due to the gauge and translation invariance of \eqref{nls}. 
\begin{definition}\label{def:stability}
  We say that the standing wave
  $u(t,{\bf x})=e^{i\omega t}\phi_\omega({\bf x})$ is orbitally stable, if
  for all $\varepsilon>0$, there exists $\delta>0$ such that if $u_0\in
  H^1$ satisfies
  \[\|u_0-\phi_\omega\|_{H^1}\le \delta,\]
  then the
  solution to \eqref{nls} with $u_{\mid t=0}=u_0$ satisfies
  \begin{equation*}
    \sup_{t\in \R}\inf_{{\theta\in \R}\atop{{\bf x}_0\in
      \R^d}}\left\|u(t, \cdot)-e^{i\theta}\phi_\omega(\cdot
      -{\bf x}_0)\right\|_{H^1}\le \varepsilon.
  \end{equation*}
  Otherwise, the standing wave is said to be unstable. 
\end{definition}
In view of Fig.~\ref{ME}, we have that $\frac{\partial}{\partial \omega}M(\phi_\omega)>0$ for the family of profiles \eqref{phi}. The well-known Vakhitov and Kolokolov stability criterion \cite{VaKo} 
therefore indicates orbital stability of these 1D solitary waves. A rigorous proof of this fact was given in \cite{Oh}, and we also refer to \cite{GSS, We} for a more general mathematical theory 
concerning the (in-)stability of standing waves.

\begin{figure}[htb!]
  \includegraphics[width=0.49\textwidth]{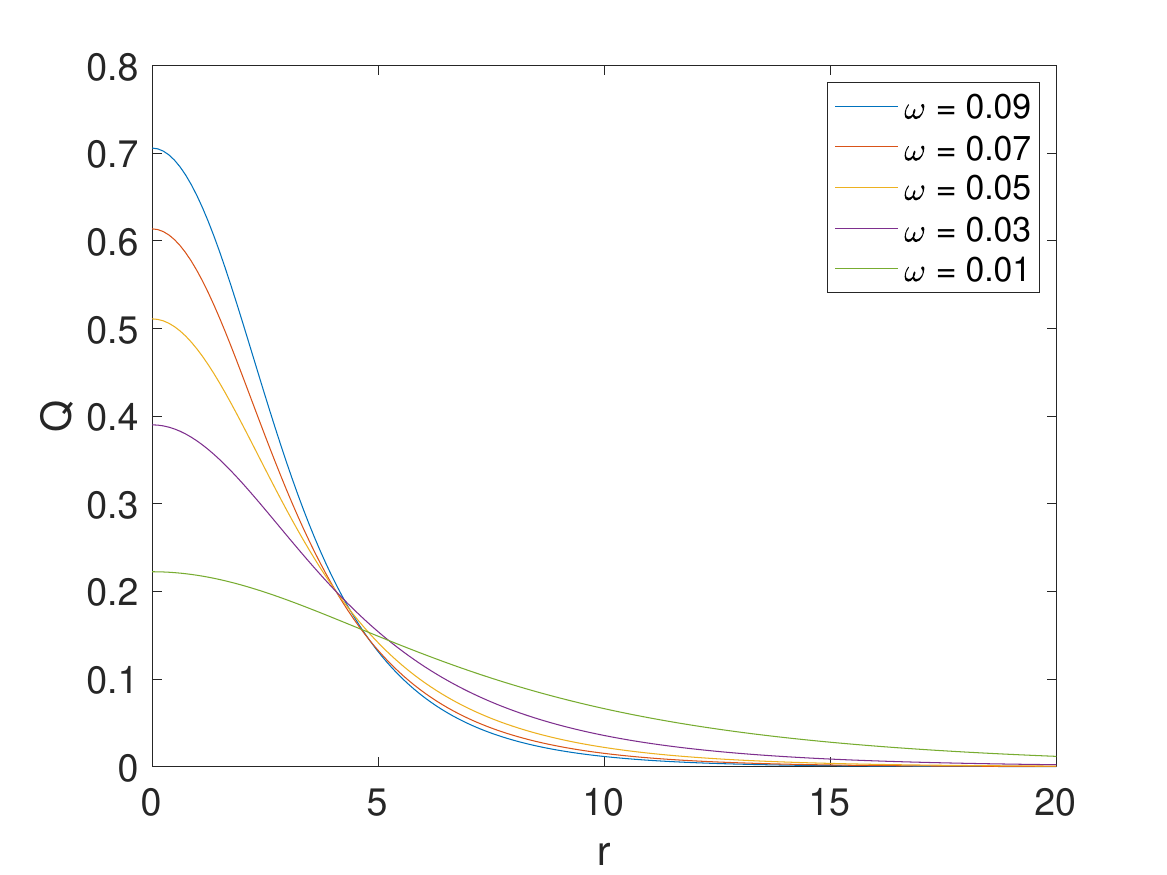}
  \includegraphics[width=0.49\textwidth]{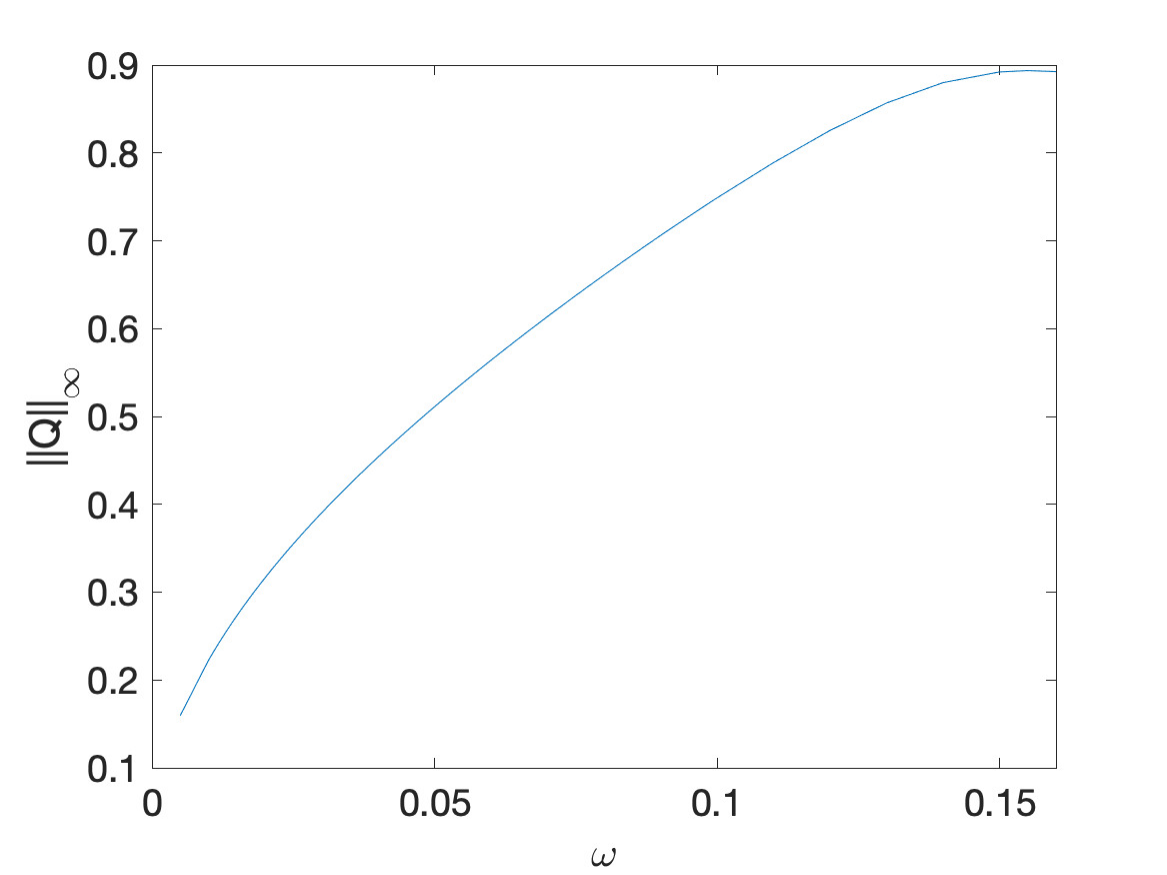}
\caption{Left: Numerically constructed 2D ground state solutions to \eqref{nls}
 for several values of $\omega$. Right: The $L^{\infty}$-norm of these 
 states as a function of $\omega$.}
 \label{NL35_d2sol}
\end{figure}

In higher dimensions, no explicit formula for profiles $\phi_\omega({\bf x})$ with ${\bf x}\in \R^d$, is known. However it can be proved that for any given $\omega  \in (0, \frac{3}{16})$, 
there exists a unique (up to gauge invariance) solution $\phi_\omega ({\bf x}) = Q_\omega(|{\bf x}|)\ge 0$, i.e., a radially symmetric, non-negative profile $Q_\omega$, which, 
together with its partial derivatives up to order two, decays exponentially to zero as $r\equiv |{\bf x}|\to \infty$,  see, e.g., \cite{CaSp} for more details. 
In 1D we clearly have (by uniqueness) that $Q_\omega$ is identical to the profile $\phi_\omega$ given in \eqref{phi}.
In 2D these (cubic-quintic) {\it ground state solutions} $Q_\omega$ have been numerically constructed in \cite{CKS}. For the convenience of the reader, we shall briefly recall 
some of its findings: 

In Fig.~\ref{NL35_d2sol} we show a plot of the 2D profiles $Q_\omega$ and their $L^\infty$-norms for several values of $\omega  \in (0, \frac{3}{16})$. 

\begin{figure}[htb!]
  \includegraphics[width=0.49\textwidth]{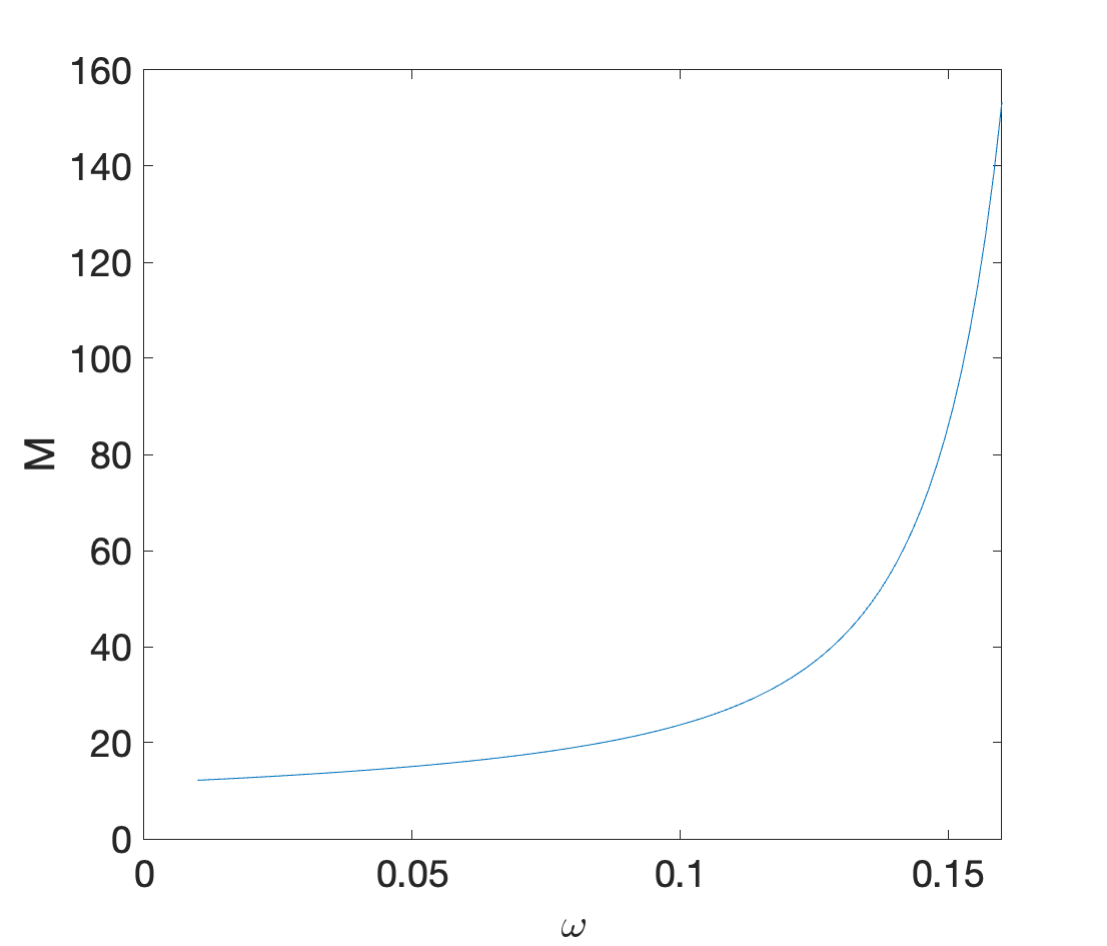}
  \includegraphics[width=0.49\textwidth]{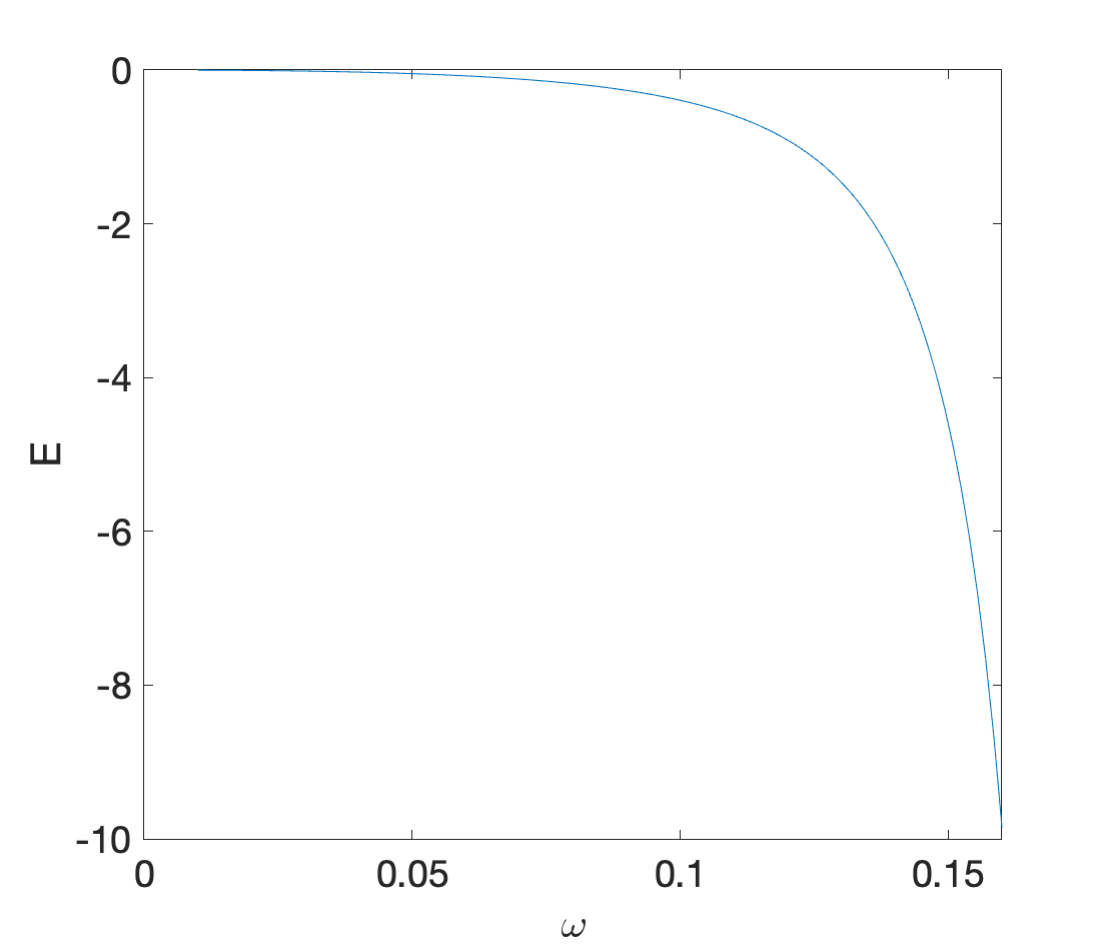}
 \caption{$M(Q_\omega)$ and $E(Q_\omega)$ as functions of $\omega$ in 2D.}
 \label{NL35_d2solmass}
\end{figure}

In Fig.~\ref{NL35_d2solmass} we show the (numerically computed) mass and the energy of these 2D solitary wave 
solutions as functions of $\omega$. Both quantities are seen to remain monotonic with $\frac{\partial}{\partial \omega}M(Q_\omega)>0$, 
just as in the 1D case above. Thus, one again expects orbital stability of the whole family of 2D ground states \cite{CKS, LeRo}. 
However, no rigorous proof of this is available to date (see \cite{CKS} for a numerical study of the long-time behavior of perturbed ground states in 2D and 3D). 

\subsection{Transverse (in-)stability of line solitary waves}\label{sec:heur}
 In this work, we are interested in an in-between situation, namely the {\it transverse (in-)stability} of the profiles \eqref{phi}. Clearly, the solution $\phi_\omega=\phi_\omega(x)$ with $x\in \R$ can be seen 
as a particular kind of 2D solution to \eqref{nls} which is constant in the transverse variable $y$, a so-called {\it line soliton}. For ${\bf x}=(x,y)\in \R^2$ such (partially constant) solutions have infinite mass and energy, which 
is why we shall instead consider \eqref{nls} on a cylindrical domain with $(x,y)\in \mathbb R_x \times \mathbb T_{L_y}$. Here $\mathbb T_{L_y} = \R/(2\pi L_y \Z)$ denotes the torus of length $L_y>0$. 
This setting is physically meaningful as it describes propagation within waveguide structures. 

The well-posedness theory for NLS on such waveguide manifolds is 
by now well understood (albeit not yet as complete as for the usual Euclidian setting), see, e.g., \cite{FYZ, Lu, TaTz, Xu} and the references therein. 
In our setting we are consequently guaranteed the existence of global in-time solutions $u\in C(\R; H^1(\mathbb R_x \times \mathbb T_{L_y}))$ 
to \eqref{nls}, depending continuously on the initial data $u_0 $. 
The definition of transverse (orbital) 
stability of line solitons \[u(t, x) = \phi_\omega(x)e^{i \omega t}, \quad \omega \in (0, \tfrac{3}{16}),\] under the evolution of a 2D nonlinear flow is thereby the same 
as given in Definition \ref{def:stability}, except that the parameter ${\bf x}_0\in \R^d$ has to be replaced by $x_0\in \R$.

Rousset and N. Tzvetkov studied these type of problems in 
a series of seminal works \cite{RoTz1, RoTz2, RoTz3} 
in which they developed a rather general mathematical framework for the transverse (in-)stability of solitary waves to dispersive equations in 2D. More
recent studies of the transverse (in-)stability of line solitons in nonlinear Schr\"odinger-type problems can be found in \cite{BIK, Ya1, Ya2}. 
In particular, it has been observed that the torus length $L_y>0$ plays 
an important role for the onset of transverse instability. To this end, one should first note that for large $L_y\gg 0$, the shape of the line soliton $\phi_\omega=\phi_\omega(x)$ is far from 
the stable, radially symmetric profile of the 2D ground state $Q_\omega=Q_\omega(r)$, while for $L\to 0$ these two become more and more similar (at least formally). In this context, $Q_\omega $ is also referred to as a 
{\it lump soliton}.

Second, the possibility of finite time blow-up for 2D cubic NLS offers a possible channel of (strong) instability of solitary waves. 
To this end, we recall that for cubic NLS in 2D, solutions with mass $M(u)> M(Q^{\rm cub}_{\omega})$ 
will blow up in finite time, cf. \cite{CazCourant}. Here, $Q^{\rm cub}_{\omega}$ is the cubic nonlinear ground state for which it is known \cite{CaSp} that its mass 
is less than the corresponding the cubic-quintic one, i.e., $M(Q^{\rm cub}_{\omega}) < M(Q_\omega)$ for all $\omega \in (0, \frac{3}{16})$. In addition, $Q_\omega^{\rm cub}$ is known 
to asymptotically yield the self-similar profile of any blow-up solution in 2D.
In our model the possibility of blow-up is excluded and hence, one might wonder if the torus length $L_y$ still plays any role. (It does, as we shall see.) 
Heuristically speaking, we expect instability of the line soliton against the formation of one or several lump solitary waves, 
if the mass 
\[
M_{\rm 2D}(\phi_\omega) \equiv \iint_{\R_x\times \mathbb T_{L_y}} \phi_\omega(x) dx\, dy = 2\pi L_y M_{\rm 1D}(\phi_\omega)= 2\pi L_y \| \phi_\omega \|_{L^2(\R_x)}^2,
\]
exceeds the ground state mass $M(Q_\omega)$ at the same value of $\omega\in (0, \frac{3}{16})$. This can of course always be achieved by choosing $L_y\ge L_{\rm crit}>0$ sufficiently large.

\begin{remark} In the papers cited above, the authors consider NLS with a single 
power-law nonlinearity. In this situation, a rescaling argument w.r.t. $y\in \R$ allows one to recast a condition on $L_y$ into an equivalent condition on $\omega$, yielding 
instability for solitary waves at frequencies $\omega > \omega_{\rm crit}>0$. 
In our case, such a rescaling argument is no longer possible. 
\end{remark}

%%%%%%%%%%%%%%%%%%%%%%%%%%%%%%%%%

\section{Numerical algorithm}\label{sec:num}

To study the NLS \eqref{nls} in two space dimensions, we will not re-apply the 
radially symmetric code used in \cite{CKS}, since we are here interested in 
a manifestly not radially symmetric situation. 

Instead we shall use the 
standard discretization of the Fast Fourier transform (FFT) in both 
$x$ and $y$. This means that in our numerics we will effectively work on a 2D torus with $x\in 
[-L_{x}\pi,L_{x}\pi]$ and $y\in [-L_{y}\pi,L_{y}\pi]$, subject to periodic boundary conditions. The parameters $L_{x}\gg L_{y}>0$ are thereby chosen large enough to avoid 
boundary effects within $u$. 
We shall consequently use $N_{x}\gg N_{y}\in \N$ FFT modes to spatially discretize the problem. 
After discretization, the NLS \eqref{nls} becomes an 
equation of the form
\begin{equation}\label{LN}
	\partial_t \widehat{u} = -i(k_{x}^2+k_{y}^2)\widehat{u}+i \mathcal F \big(|u|^4u-|u|^2u\big),
\end{equation}
where $\widehat{u}\equiv \mathcal F(u)$ denotes the FFT of $u$ and 
$L_{x}k_{x}, L_{y}k_{y}\in \mathbb Z$ are the 
Fourier variables dual to $x$ and $y$, respectively. 

There are many efficient integration schemes for equations of the form \eqref{LN}, see, e.g. \cite{etna} for a recent comparison. 
In the present paper we shall apply the composite fourth order Runge-Kutta scheme due to Driscoll \cite{Dr}. 
The quality of the code is thereby controlled as in \cite{etna}, i.e., using the 
relative conservation of the numerically computed energy 
\[
\Delta_E(t)= \left|\frac{E(t)}{E(0)}-1\right| ,
\]
together with the observed decrease of the FFT coefficients w.r.t. the index $n\in [0, N_{x,y}]$. Due to unavoidable numerical errors, the 
energy $E(t)$ will not remain conserved by the time-evolution of \eqref{LN}, and thus, the evolution of $\Delta_E$ serves a basic control parameter for the accuracy of the code.

To test the code, we first discretize
the line solitary wave $\phi_\omega$ at 
$\omega=0.1$ with $N_{x}=2^{10}$ FFT modes for $L_{x}=40$, and 
$N_{y}=2^{5}$ FFT modes for $L_{y}=3$. We then let it propagate in time using $N_{t}=10^{3}$ time steps for 
$t\in [0,1]$. The relative conservation of the energy 
at $t=1$ is observed to be $\Delta_E=\mathcal O(10^{-15})$. The observed $L^\infty$-difference between the 
numerical solution at $t=1$ and the exact solution at the final time 
\[u(t,x)_{\mid t=1} = \phi_{\omega = 0.1 }(x) e^{0.1 i }\] is of the 
order of $10^{-16}$, i.e., of the order of machine precision. 

Whereas this tests the performance of the code for the evolution of line solitary waves, it does not control the 
$y$-dependence of possible perturbations. Hence, as a second test, we consider initial data 
$u_0 = Q_{\omega =0.1}$, the numerically constructed ground state in 2D. Due to the radial 
symmetry of the solution $u$ in this case, we take $L_{x}=L_{y}=10$ and discretize this domain using $N_{x}=N_{y}=2^{8}$ Fourier modes. 
We then let these particular initial data propagate using $10^{4}$ time steps for 
$t\in [0,1]$. The $L^\infty$-difference between the 
numerical solution $u$ at $t=1$ and 
\[u(t,x)_{\mid t=1} =Q_{\omega=0.1}(|{\bf x}|) e^{0.1i}\] 
 is of the order of $10^{-15}$, i.e., again of the order of machine precision 
(the relative energy conservation $\Delta_E$ is of the same order). 

\begin{remark} Note that this also tests the quality of the numerically constructed ground state solution $Q_\omega$ which 
was obtained iteratively in \cite{CKS} with a residual of the order of 
$10^{-12}$. 
\end{remark}

In summary, these tests show that our numerical algorithm is able to solve the 
cubic-quintic NLS in 2D to the order of machine precision, 
and that energy conservation provides a valid control of the 
resolution in time. 

%%%%%%%%%%%%%%%%%%%%%%%%%%%%%%%%%%%%%%%%%%

\section{Blow-up of perturbed line solitons in the cubic NLS}\label{sec:cubic}

Since we expect the transverse instability of line solitons in cubic-quintic NLS to be driven by the 
underlying blow-up dynamics in the purely cubic case, we shall briefly study the latter situation in this section. 

The focusing cubic NLS is given by
\begin{equation}\label{cubic}
\left\{\begin{aligned}
	  i\partial_t u +\Delta u & = - |u|^2u ,\\
   u_{\mid
  t=0} & = u_0,
\end{aligned}
\right.
\end{equation}
and we can expect finite time blow-up of its solution $u$ whenever $M(u_0)>M(Q^{\rm cub}_\omega)$. 
Line solitary waves for cubic NLS correspond to $y$-independent profiles 
\[
\phi^{\rm cub}_{\omega}(x) = \sqrt{2\omega}\, \mbox{sech}(\sqrt{\omega}x), \quad \omega \in (0, \infty).
\]
The corresponding $L^\infty$ norm is $\| \phi^{\rm cub}_\omega\|_{L^\infty}= \sqrt{2\omega}$ and we also note that due to 
the $L^2$-scaling invariance of the problem any positive frequency $\omega$ is admissible in this case. In particular, the ground state mass $M(Q_\omega^{\rm cub})$ is the same for 
all values of $\omega$. 

In the following, we shall consider initial data of the form 
\begin{equation}
	u_{0, \pm}(x,y) = \phi^{\rm cub}_{\omega}(x) \pm \frac{ \sqrt{2\omega}}{10} e^{-(x^{2}+y^{2})}
	\label{pertcub},
\end{equation}
i.e., transverse perturbations of $\phi_\omega$ with a size of the order of $10\%$ of $\| \phi^{\rm cub}_\omega\|_{L^\infty}$. 
Numerically, we will work with the parameters $L_{x}=100$, $L_{y}=2$ and $N_{x}=2^{12}$, 
$N_{y}=2^{7}$ FFT modes. We also take $N_{t}=10^{4}$ time steps for $t\in [0, 100]$. 

\begin{figure}[htb!]
  \includegraphics[width=0.49\textwidth]{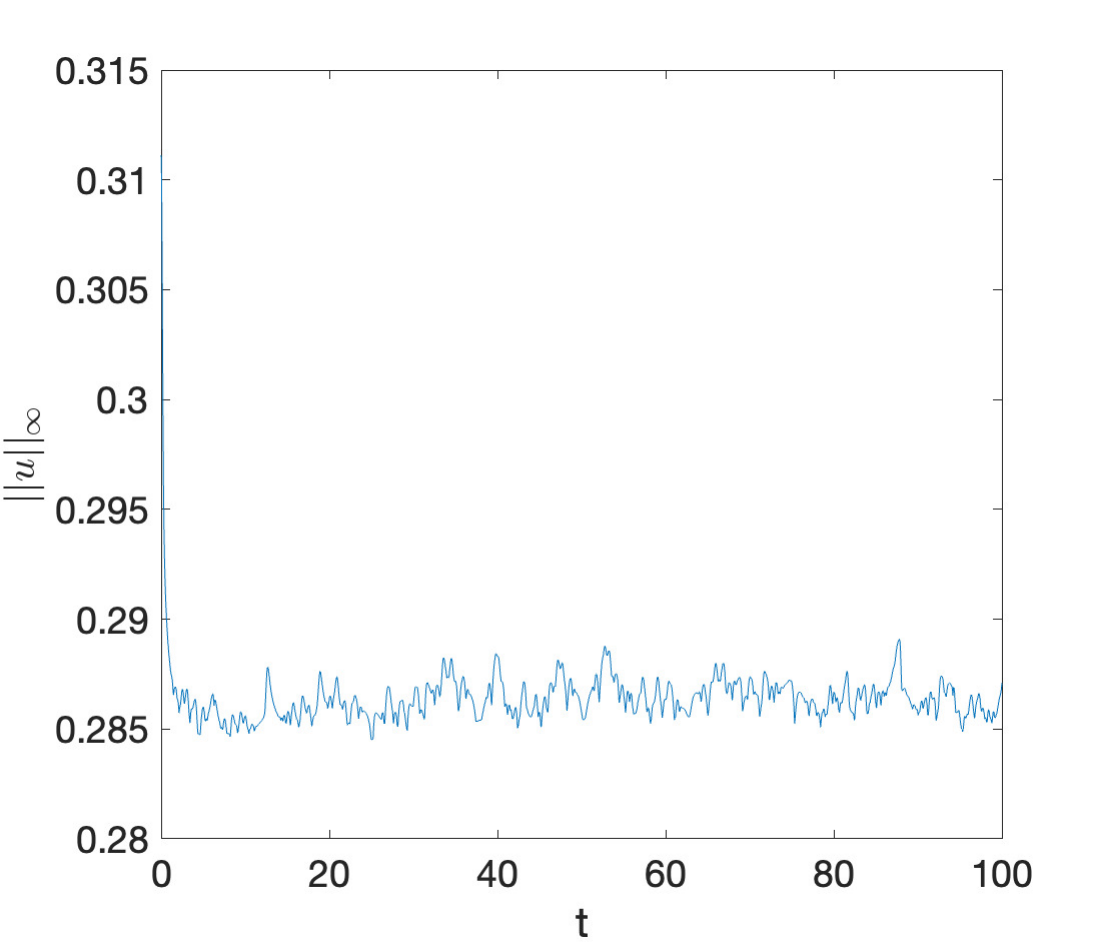}
  \includegraphics[width=0.49\textwidth]{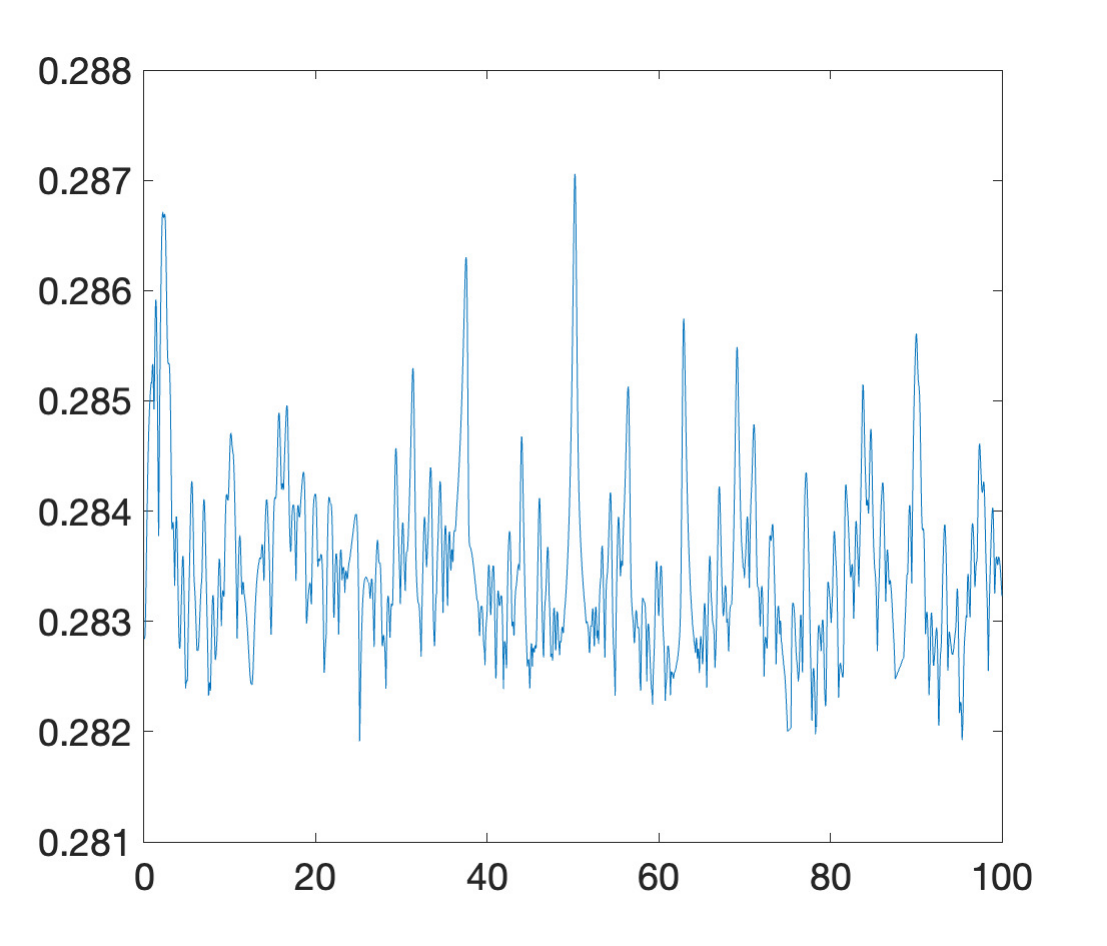}
 \caption{Time-evolution of the $L^{\infty}$ norm of the solution $u$ 
 to (\ref{cubic}) for initial 
 data \eqref{pertcub} with $\omega=0.04$ (on the left for $u_{0,+}$ and on the right for $u_{0,-}$). }
  \label{pertcub004}
\end{figure}

First we shall consider the case $\omega=0.04$, for which we find $M(u_{0, +}) \approx 10.10$ and 
$M(u_{0, -}) \approx 10.05$, respectively. Since both of these are smaller than the ground state mass
$M(Q^{\rm cub}_{\omega=0.1})\approx 11.7$ the solution $u$ to \eqref{cubic} 
will not blow-up. Numerically this is confirmed by the time-evolution of the $L^{\infty}$ norm of $u$ which is shown in Fig.~\ref{pertcub004} for both choices of the $\pm$ sign. 
In addition, the perturbed line soliton appears to remain stable on the considered time-interval $t\in [0, 100]$. Since we are working on a (large) torus $\mathbb T_{L_x}\times \mathbb T_{L_y}$, 
radiation cannot escape to infinity and hence, no final state is reached as can be observed from the persistent oscillations in Fig.~\ref{pertcub004}. 

\begin{figure}[htb!]
  \includegraphics[width=0.49\textwidth]{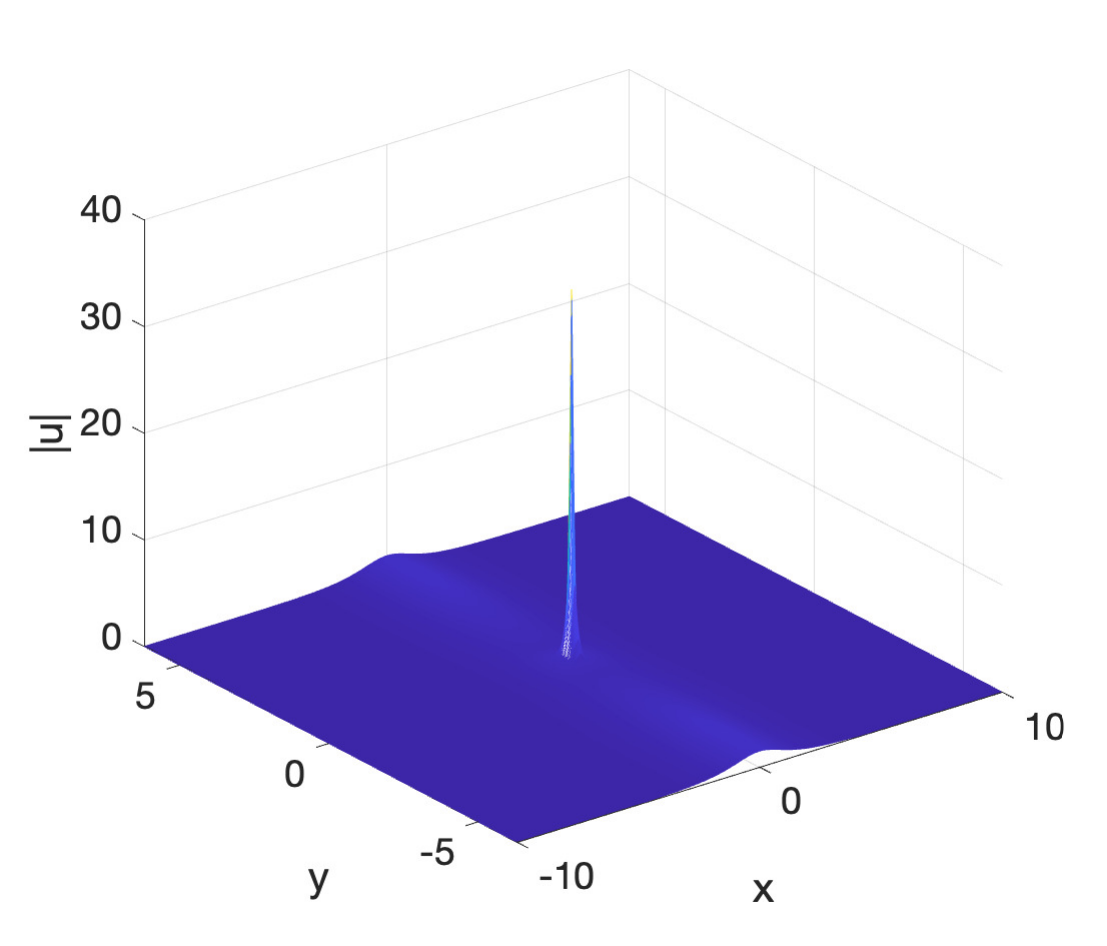}
  \includegraphics[width=0.49\textwidth]{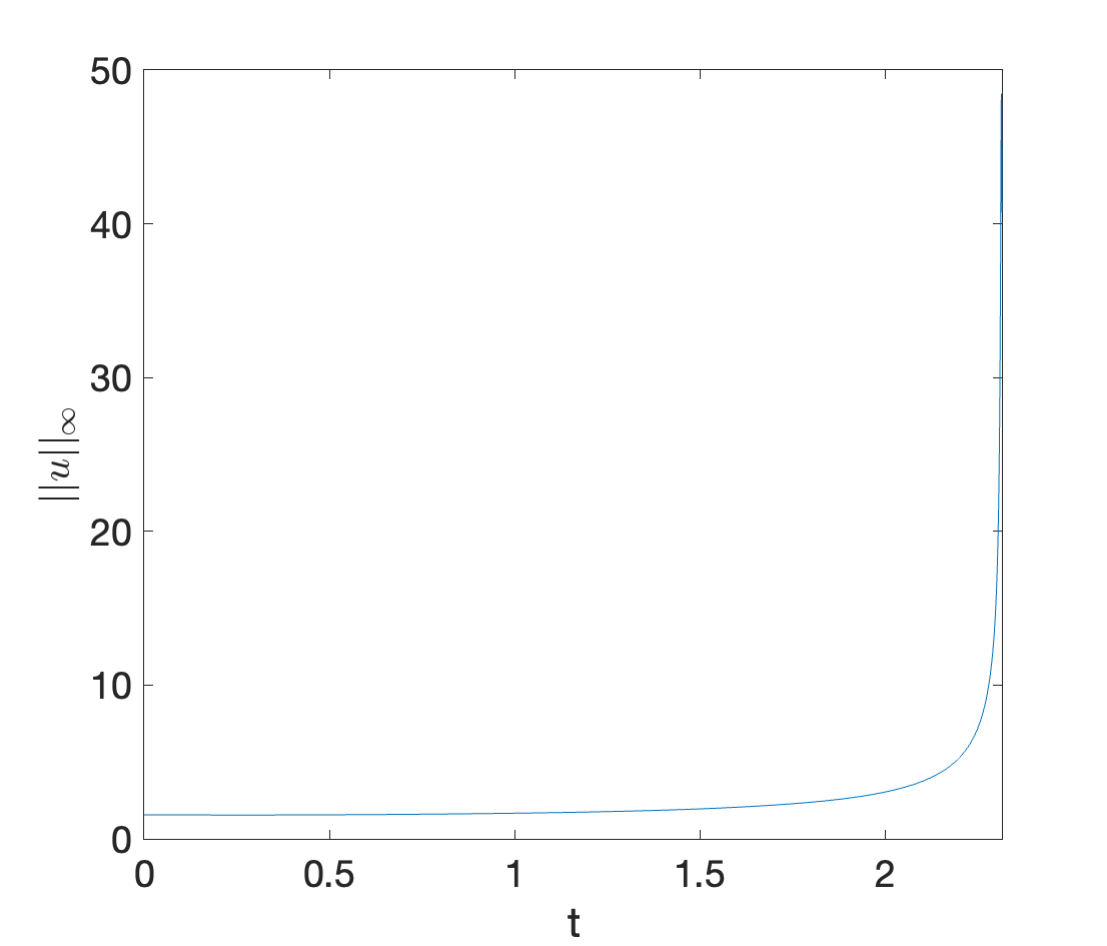}
 \caption{Left: Solution $|u|$ to the cubic NLS at time $t=2.316$ corresponding to initial 
 data $u_{0, +}$ with $\omega=1$. Right: Time-evolution of the $L^{\infty}$ norm of the 
 solution $u$ on the $y$-axis.}
  \label{pertcubp}
\end{figure}

If we consider initial data of the form \eqref{pertcub} with $\omega=1$, we find that $M(u_{0, +}) \approx 51.35$  
and $M(u_{0, -}) \approx 49.25$. In both cases the mass of the initial data $u_0$ is more than four times larger than $M(Q_{\omega =1}^{\rm cub})\approx 11.7$ which will cause the solution $u$ to blow-up. 
Indeed in the case with initial data $u_{0,+}$ the initial perturbation grows 
monotonically in time to large values $\|u (t, \cdot, \cdot)\|_{L^\infty} \approx 50$. The code is stopped at $t=2.316$, where the 
relative conservation of the energy $\Delta_E(t)$ drops below $10^{-3}$ and hence, too much accuracy is lost. In order to investigate the blow-up in more detail (which is beyond the scope of this paper), a 
much higher numerical resolutions would need to be used. 
We show the modulus of the solution at the final recorded 
time in Fig.~\ref{pertcubp} on the left. The time-evolution of the $L^{\infty}$ norm of $u$ on the right of the same figure indicates a blow-up. 

\begin{figure}[htb!]
  \includegraphics[width=0.49\textwidth]{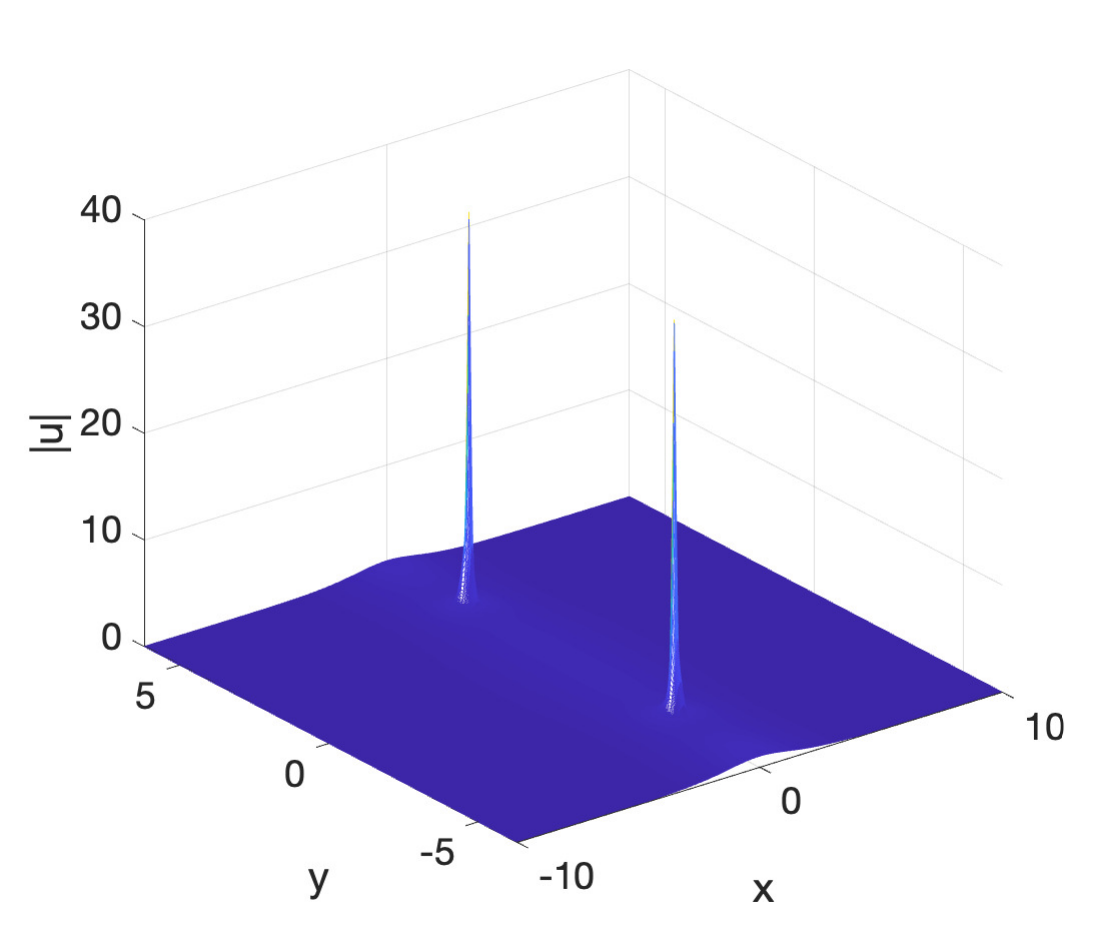}
  \includegraphics[width=0.49\textwidth]{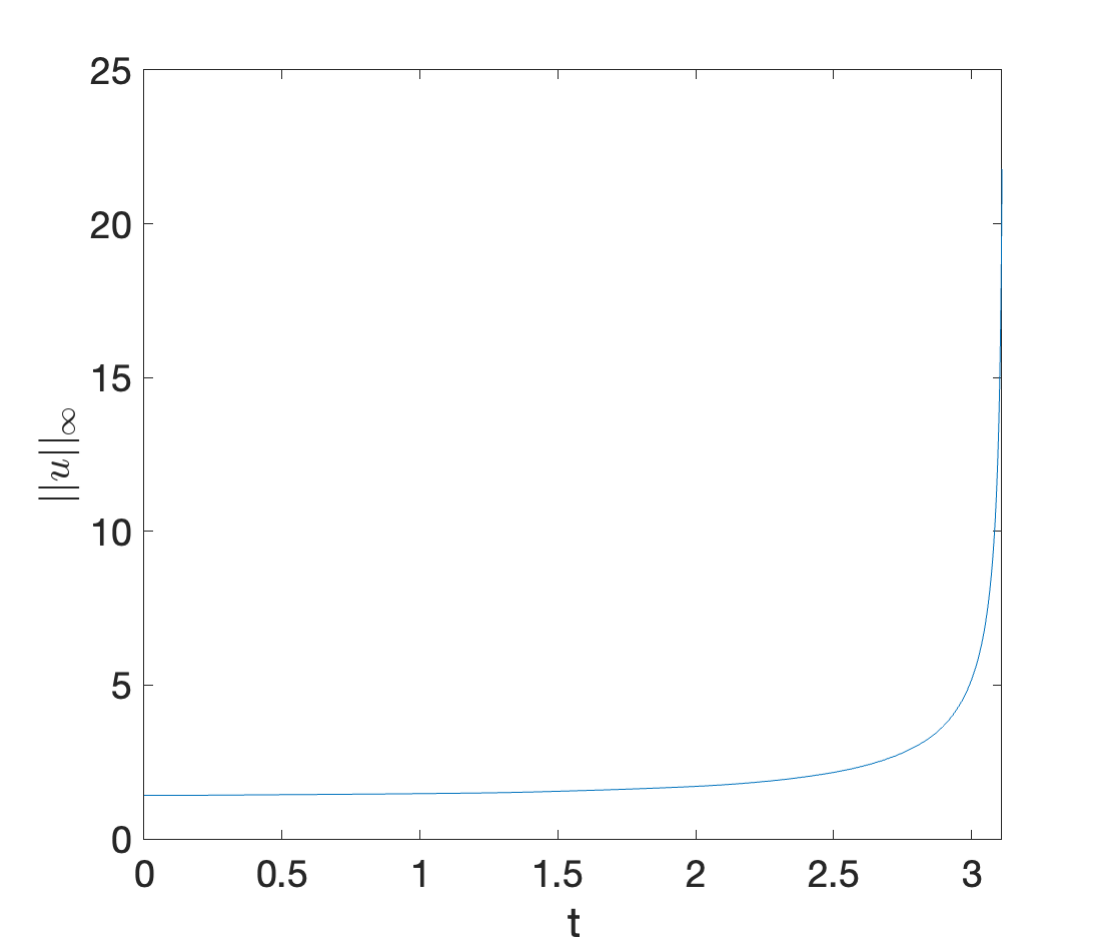}
 \caption{Left: Solution $|u|$ to the cubic NLS at time $t=3.12$ corresponding to initial 
 data $u_{0, -}$ with $\omega=1$. Right: Time-evolution of the $L^{\infty}$ norm of the 
 solution $u$ on the $y$-axis.}
  \label{pertcubm}
\end{figure}

The situation for initial data $u_{0,-}$ is depicted in the Fig.~\ref{pertcubm}. 
In this case the code stops slightly later at $t=3.12$ and the solution at the final 
time shows blow-up at two well-separated peaks (formed from 
the rims of the initial hole). The $L^{\infty}$ norm of the solution 
as a function of time is shown on the right of the same figure. 

%%%%%%%%%%%%%%%%%%%%%%%%%%%%%%%%%%%

\section{Cubic-quintic NLS: the stable regime}\label{sec:stab}

We consider transverse perturbations of cubic-quintic line solitons given by initial data of the form
\begin{equation}\label{pert}
	u_0(x,y) = \phi_{\omega}(x)+\lambda e^{-(x^{2}+y^{2})},
\end{equation}
where $\lambda\in \mathbb{R}$ is chosen in a way to ensure a perturbation of the order of $10\%$ when compared to the $L^\infty$ norm of $\phi_\omega$. 
The torus parameters $L_{x}$, $L_{y}$ are taken large enough to guarantee that the FFT coefficients of $u_0$ decrease exponentially and the perturbation is numerically equal to zero at the boundary.

As a first example we take $\omega=0.1$ and $\lambda= 0.05$. 
The numerical parameters are $L_{x}=40$, 
$L_{y}=2$, with $N_{x}=2^{10}$, $N_{y}=2^{7}$ and $N_{t}=10^{3}$ 
time-steps for $t\in [0,20]$. 
In this case $M_{\rm 2D}(\phi_{\omega=0.1}) = 4\pi M_{\rm 
1D}(\phi_{\omega=0.1}) \approx 20.22$, while 
the mass of the nonlinear ground state is $$M(Q_{\omega =0.1})\approx 23.74.$$
Thus we are in a regime, where $M_{\rm 2D}(\phi_{\omega=0.1})< M(Q_{\omega=0.1})$ for which our heuristic argument given in Section \ref{sec:heur} predicts stability.

\begin{figure}[htb!] \label{figom01p}
  \includegraphics[width=0.7\textwidth]{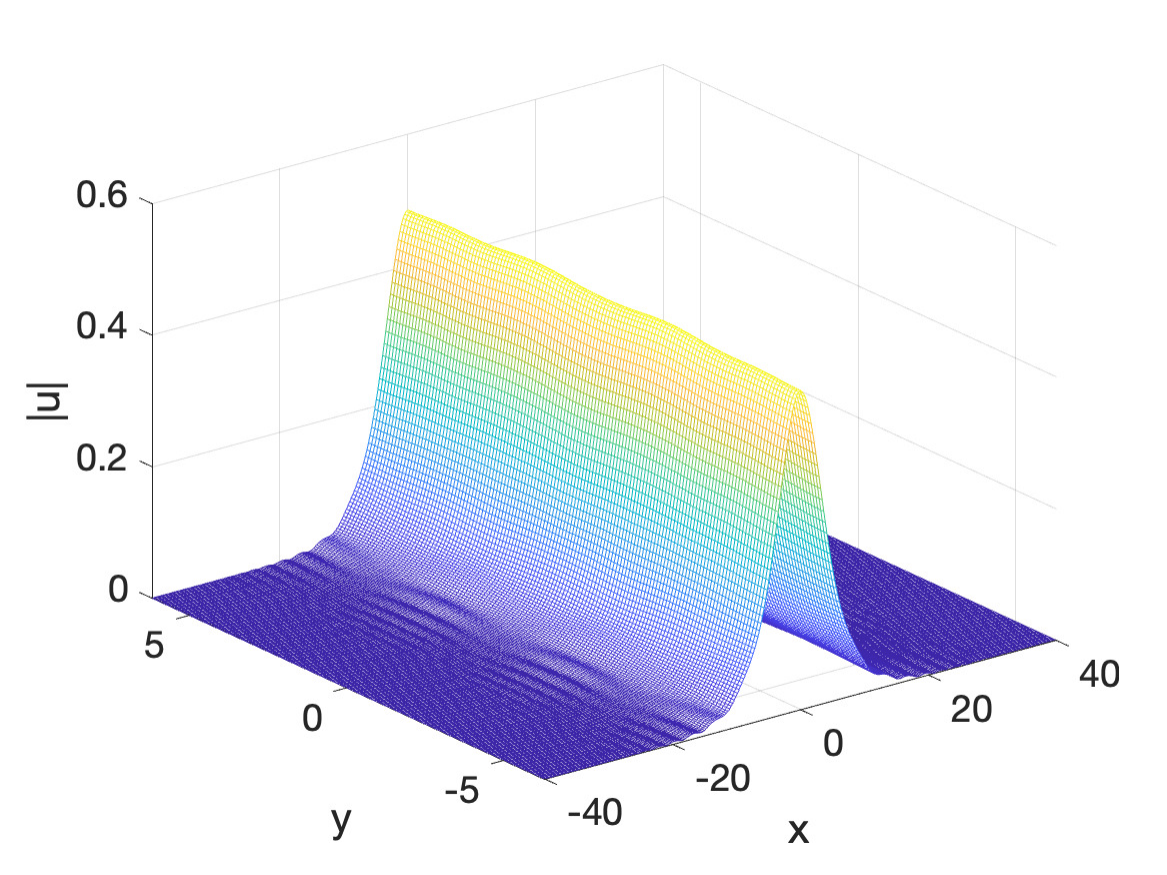}
 \caption{Solution  $|u|$ of (\ref{nls}) at time $t=5$ for initial 
 data (\ref{pert}) with $\omega=0.1$, $\lambda= 0.05$ and $L_y =2$. }
\end{figure}

\begin{figure}[htb!]
  \includegraphics[width=0.49\textwidth]{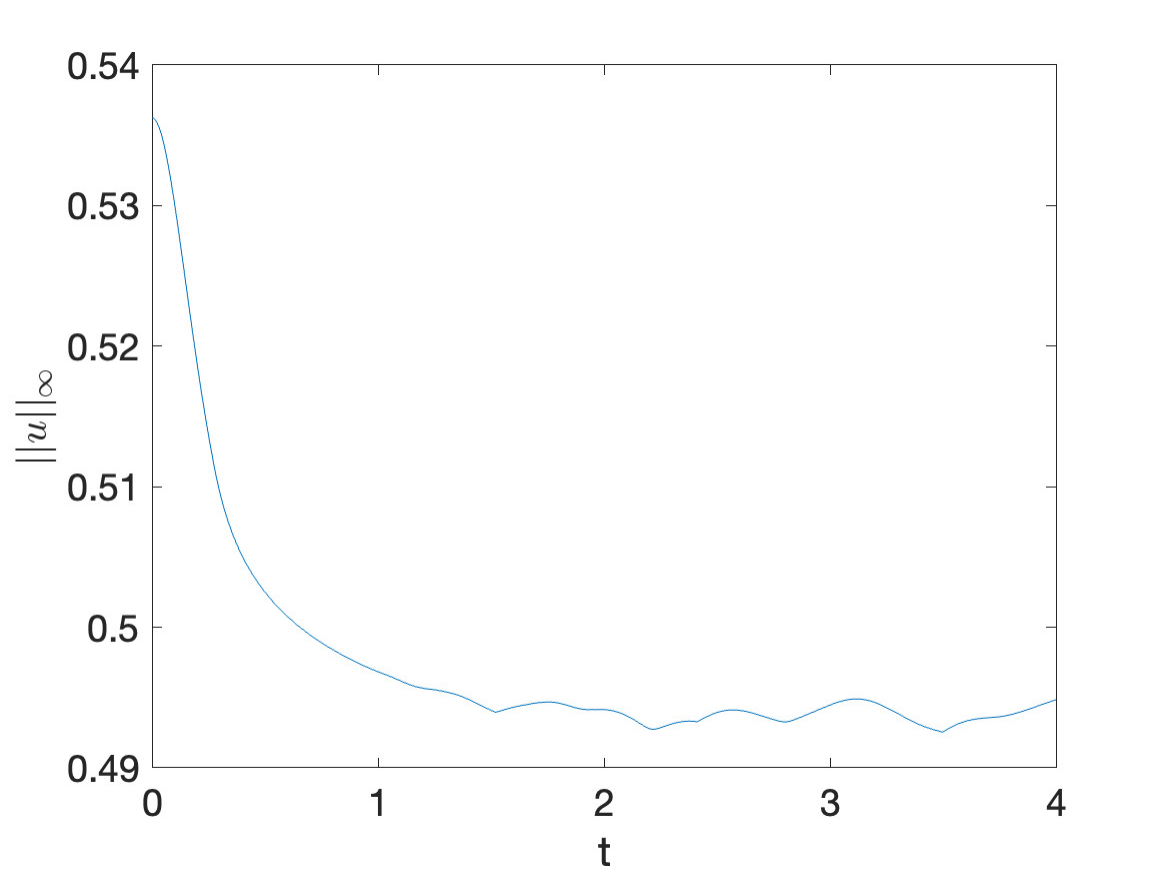}
  \includegraphics[width=0.49\textwidth]{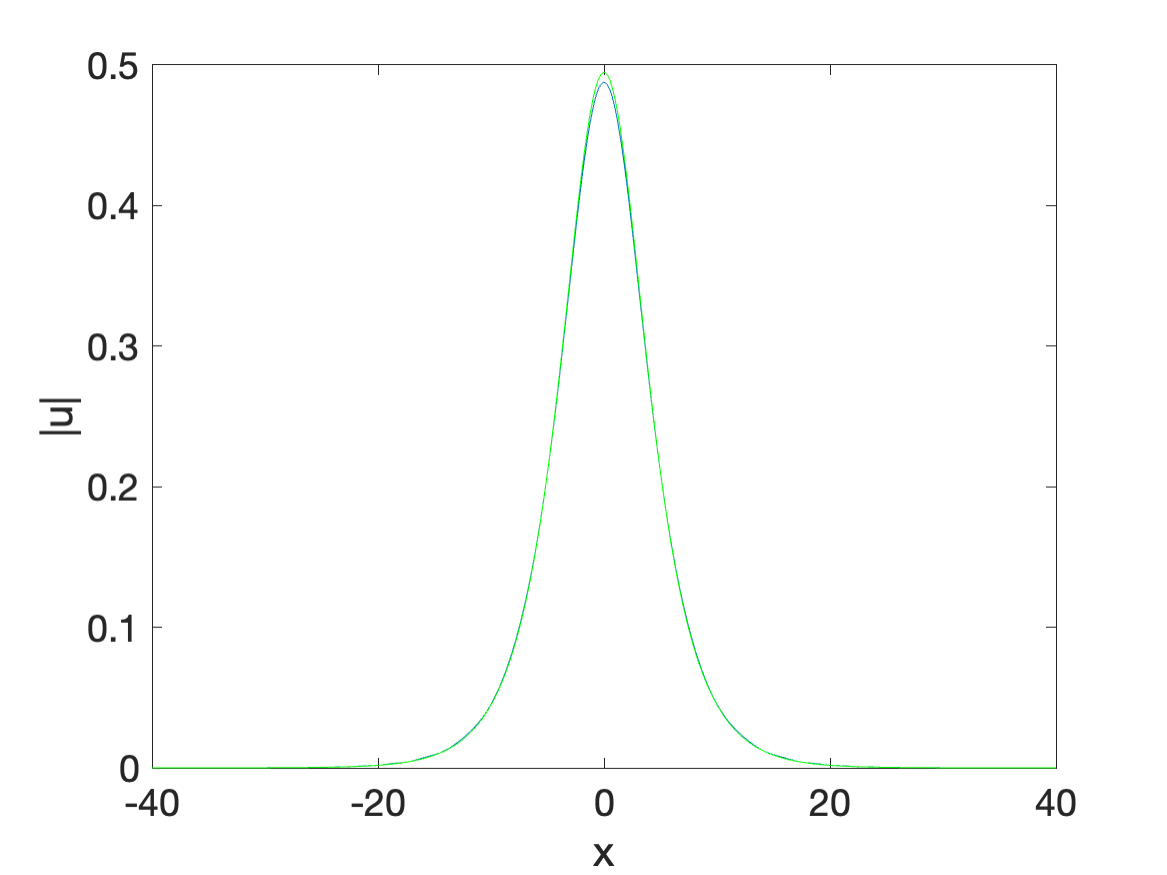}
 \caption{Left: Time-evolution of the $L^{\infty}$ norm of the solution $u$ to (\ref{nls}) for initial 
 data \eqref{pert} with $\omega=0.1$, $\lambda=0.05$ and $L_y =2$. 
 Right: The modulus of $u$ at time $t=5$ on the $y$-axis in blue 
 together with a fitted solitary wave $\phi_{\omega_\ast}$ in green. }
  \label{figom01pinf}
\end{figure}

The modulus of the solution $u$
at time $t=5$ is shown in Fig.~\ref{figom01p}, while the (short) time evolution of $\| u(t, \cdot, \cdot)\|_{L^\infty}$ is shown on the left of Fig.~\ref{figom01pinf}.

The solution clearly appears to retain its shape as a line solitary wave $\phi_\omega$ with some new, unknown 
value $\omega=\omega_\ast $. This interpretation is confirmed by fitting (a $y$-averaged value of) the $L^\infty$ norm 
of $u(t, x,y)$ at $t=5$ to the $L^\infty$ norm of the family of explicit profiles $\phi_\omega$ given in \eqref{phi}. The latter is seen to be equal to
\[
\| \phi_\omega \|_{L^\infty} = \phi_\omega (0)  = \sqrt{3} \left(\tfrac{1}{4}-\sqrt{\tfrac{1}{16}-\tfrac{\omega}{3}}\right)^{1/2}.
\]
Doing so, we numerically identify $\omega_\ast \approx 0.1017$, i.e., 
a value slightly large than for the unperturbed profile $\phi_{\omega =0.1}$. We show the solution on 
the $y$-axis in blue together with the fitted solitary wave in green 
on the right of Fig.~\ref{figom01pinf}. It can be seen that the agreement is 
excellent even though the final state in the time-evolution is not yet reached due to the 
presence of radiation. 

Further numerical simulations show that the main profile within $u(t, \cdot, \cdot) $ does not change much even if one computes for considerably longer 
times $t\in [0, 100]$, see the following figure: 
\begin{figure}[htb!]
  \includegraphics[width=0.7\textwidth]{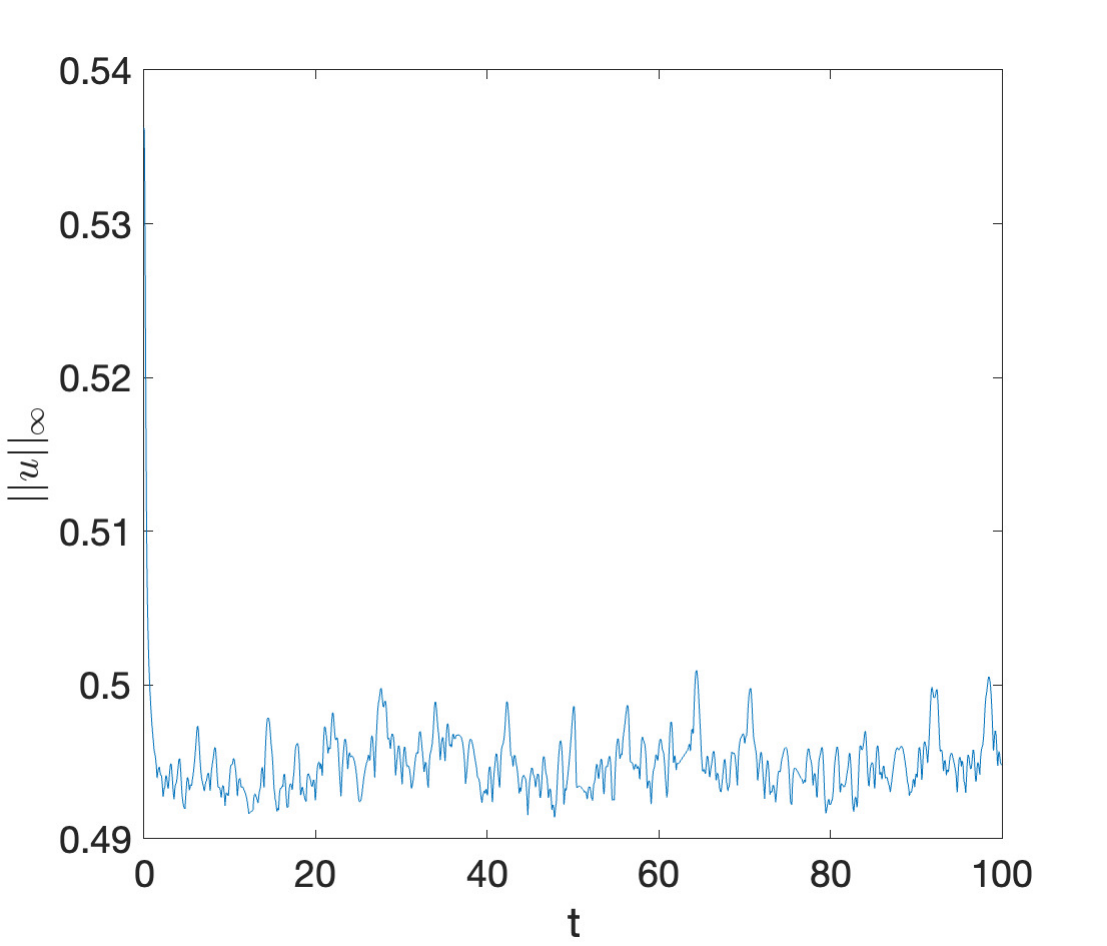}
  \caption{Long-time evolution of the $L^{\infty}$ norm of the solution $u$ to (\ref{nls}) for initial 
 data \eqref{pert} with $\omega=0.1$, $\lambda=0.05$ and $L_y =2$. }
\end{figure}

\begin{figure}[htb!]
  \includegraphics[width=0.7\textwidth]{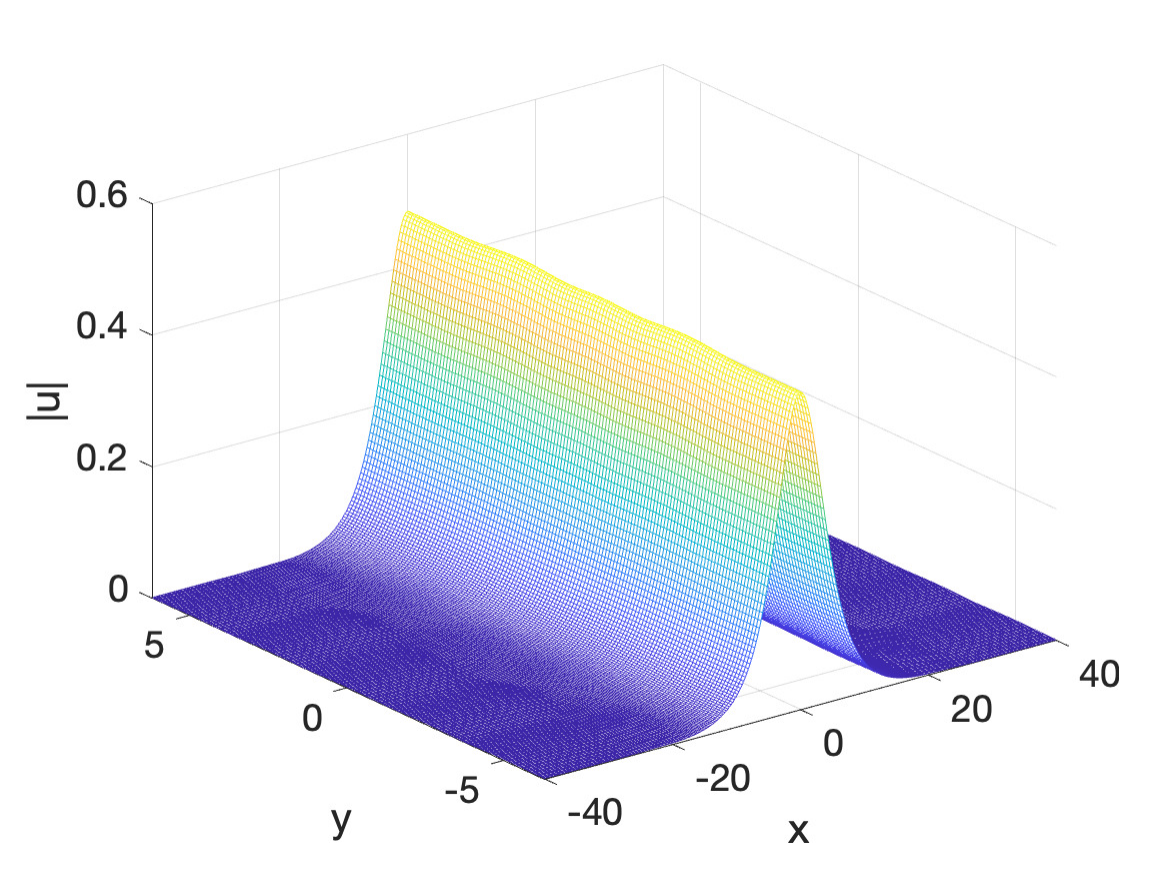}
 \caption{Solution  $|u(t,x,y)|$ to (\ref{nls}) at time $t=20$ for initial 
 data (\ref{pert}) with $\omega=0.1$, $\lambda= - 0.05$ and $L_y =2$. }
 \label{figom01m}
\end{figure}

The situation remains similar for initial data \eqref{pert} with $\lambda = -0.05$, i.e., for a 
perturbations with even smaller mass, keeping us firmly within the stable regime. We 
use the same numerical parameters as before and compute up to time $t = 20$. The modulus of the solution at the 
final time is shown in Fig.~\ref{figom01m}. Once more it clearly 
correspond to a line solitary wave plus radiation (which is seen to be markedly smaller than before).

\begin{figure}[htb!]
  \includegraphics[width=0.49\textwidth]{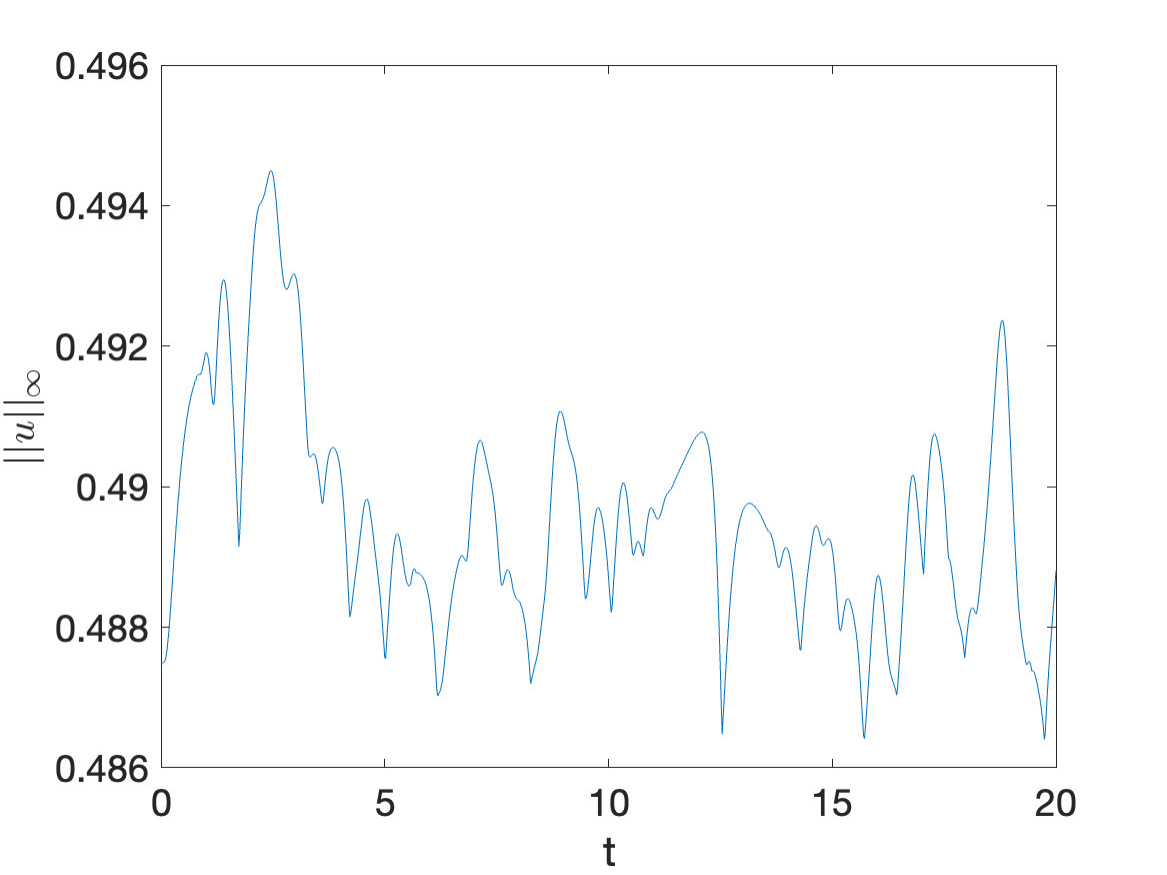}
  \includegraphics[width=0.49\textwidth]{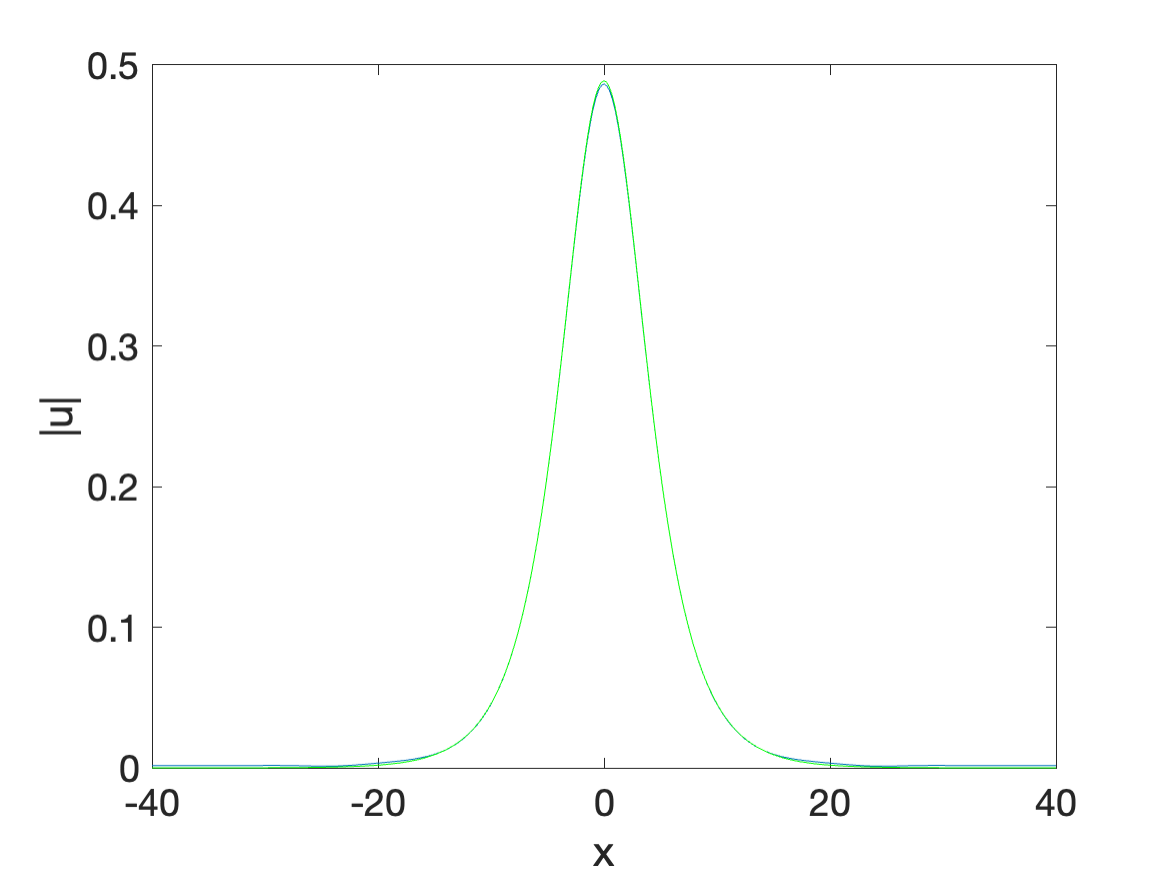}
 \caption{Left: Time-evolution of the $L^{\infty}$ norm of the solution $u$ to (\ref{nls}) for initial 
 data \eqref{pert} with $\omega=0.1$, $\lambda= - 0.05$ and $L_y =2$. Right: The modulus of $u$ at time $t=20$ on the $y$-axis in blue 
 together with a fitted solitary wave $\phi_{\omega_\ast}$ in green. }
  \label{figom01minf}
\end{figure}

As before, this interpretation is confirmed by fitting the $L^{\infty}$ norm of 
the solution $u$ at $t=20$ to profile \eqref{phi}, see the left of Fig.~\ref{figom01minf}. We thereby find the new value 
$\omega_\ast =0.1004$. The solution on the $y$-axis for $t=20$ is shown 
together with the fitted solitary wave in green, both curves being in 
excellent agreement. 

Once more the situation does not change much if we compute for even  
longer times than those shown in the figures above.

%%%%%%%%%%%%%%%%%%%%%%%%%%%%%%%%%%%%%%%%

\section{Cubic-quintic NLS: the unstable regime}\label{sec:unstab} 

To see whether the length of the period in the $y$-direction plays any role, we consider 
the same situation as in Section \ref{sec:stab}, i.e., perturbations of $\phi_\omega$ by 2D Gaussian functions, but with a larger value of $L_{y}$. 

\begin{figure}[htb!]
  \includegraphics[width=0.7\textwidth]{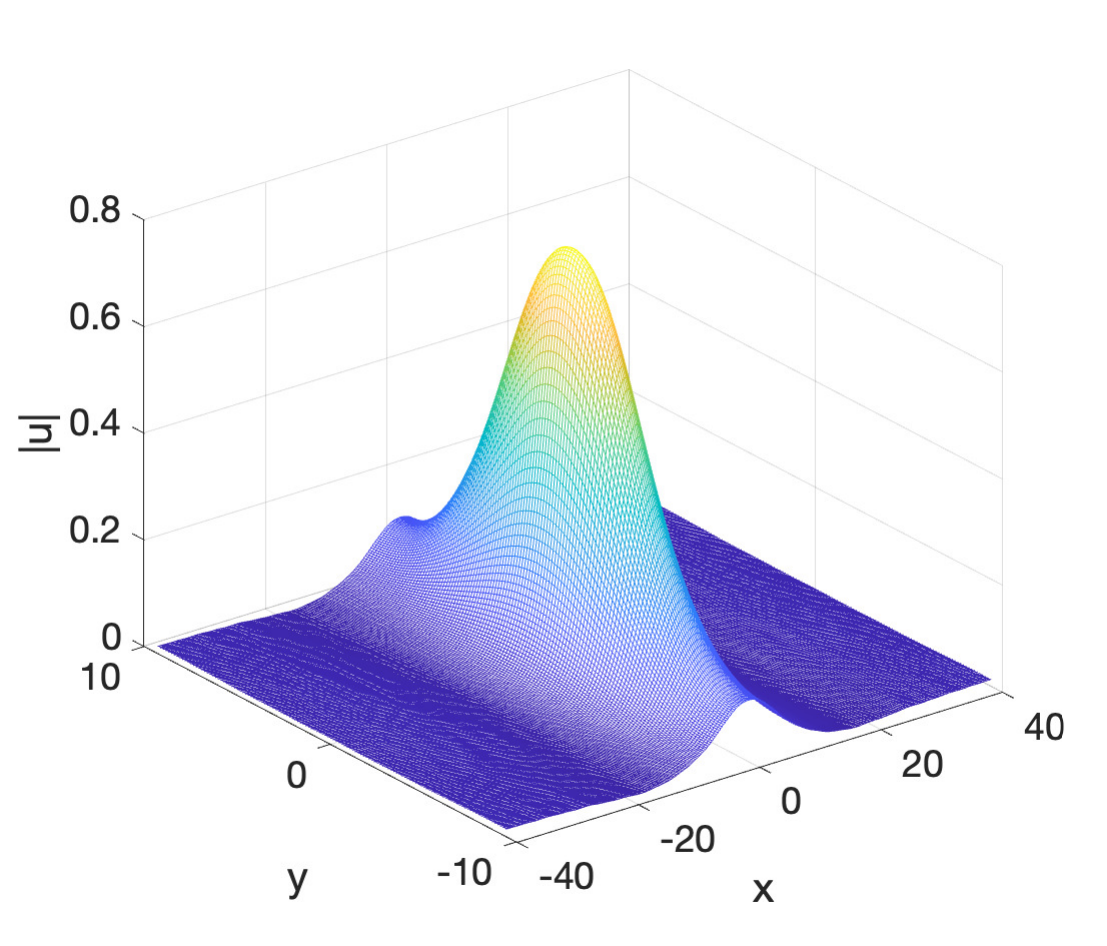}
 \caption{Solution $|u|$ of (\ref{nls}) at time $t=500$ for initial 
 data (\ref{pert}) with $\omega=0.1$, $\lambda=0.05$ and $L_y =3$. }
 \label{figom01_3p}
\end{figure}

Concretely, we use $L_{x}=150$, 
$L_{y}=3$, with $N_{x}=2^{12}$, $N_{y}=2^{7}$ and $N_{t}=5\cdot 10^{4}$ 
time-steps for $t\in [0,500]$. In this case  
\[M_{\rm 2D}(\phi_{\omega = 0.1}) = 6\pi M_{\rm 1D} 
(\phi_{\omega =0.1}) \approx 30.33> M(Q_{\omega =0.1})\approx 23.74.
\] 
We are thus in a regime where we expect instability. Indeed, for $\lambda = 0.05$ the unstable nature of $\phi_\omega$ is clearly visible in Fig.~\ref{figom01_3p}, where the modulus of the solution at time $t=500$ is shown. 

\begin{figure}[htb!]
  \includegraphics[width=0.49\textwidth]{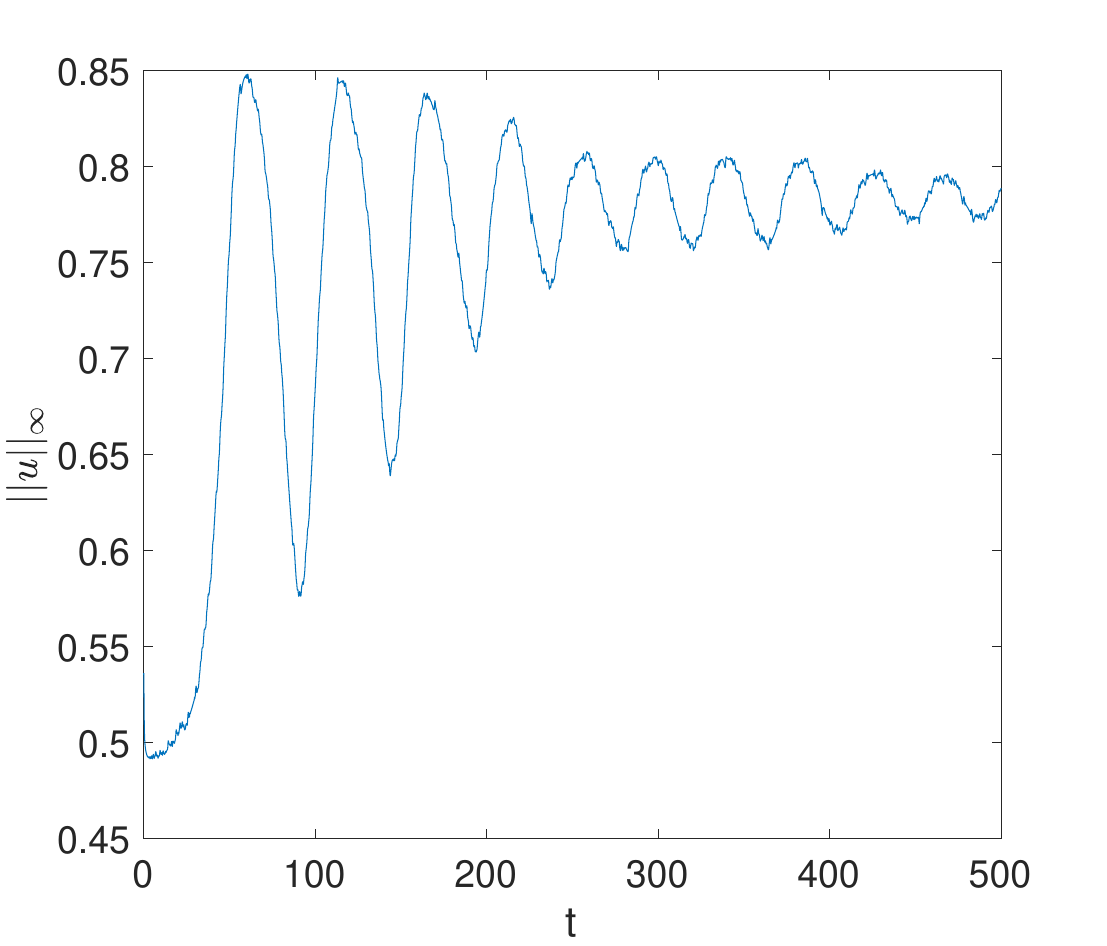}
  \includegraphics[width=0.49\textwidth]{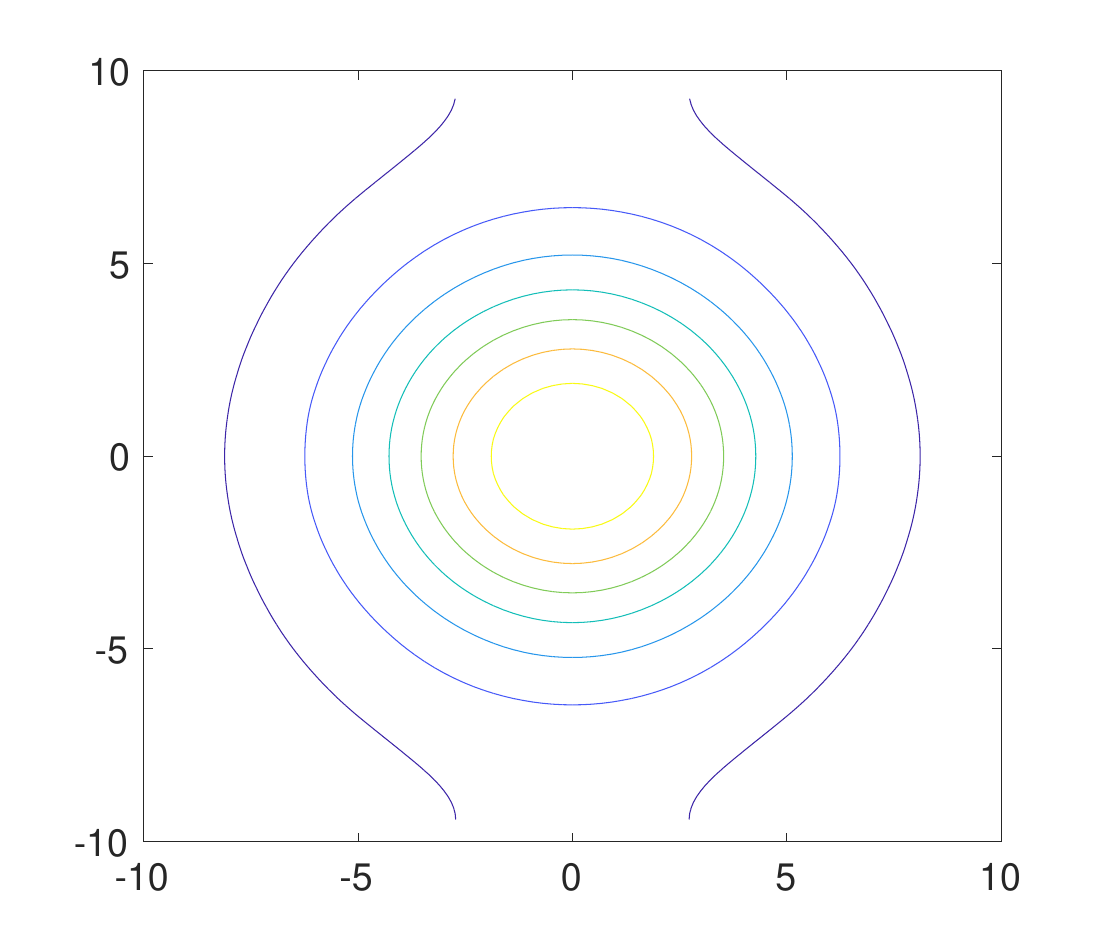}
 \caption{Left: Time-evolution of the $L^{\infty}$ norm of the solution to (\ref{nls}) for the initial 
 data (\ref{pert}) with $\omega=0.1$, $\lambda=0.05$ and $L_y =3$. Right: 
 Contour plot of the modulus of the solution $u$ at the final 
 time $t=500$. }
 \label{figom01_3pinf}
\end{figure}

The $L^{\infty}$ norm of the solution, shown on the left of 
Fig.~\ref{figom01_3pinf}, is seen to oscillate in time with a very 
large, unknown period. The amplitude, however, is clearly decreasing and thus, we expect the solution $u$ to converge as $t\to +\infty$ to some (stable) final state, plus small radiation.

We conjecture this final state to be a lump solitary wave $Q_{\omega}$ at some yet unknown frequency $\omega =\omega_\ast$. 
This is also confirmed by the contour plot on the right of the same figure, where the final state appears to be (almost) radially symmetric. 
Since there is no explicit form of $Q_\omega$, we do not present a fit of the $L^{\infty}$ norm of 
our final solution with a $\| Q_{\omega_\ast} \|_{L^\infty} $. Note however, that Fig.~\ref{NL35_d2sol} shows $\| Q_{\omega = 0.1} \|_{L^\infty} \approx 0.75$, which is 
of the same order of magnitude as $\| u(t, \cdot, \cdot)\|_{L^\infty}$ for $t=500$.

\begin{figure}[htb!]
  \includegraphics[width=0.7\textwidth]{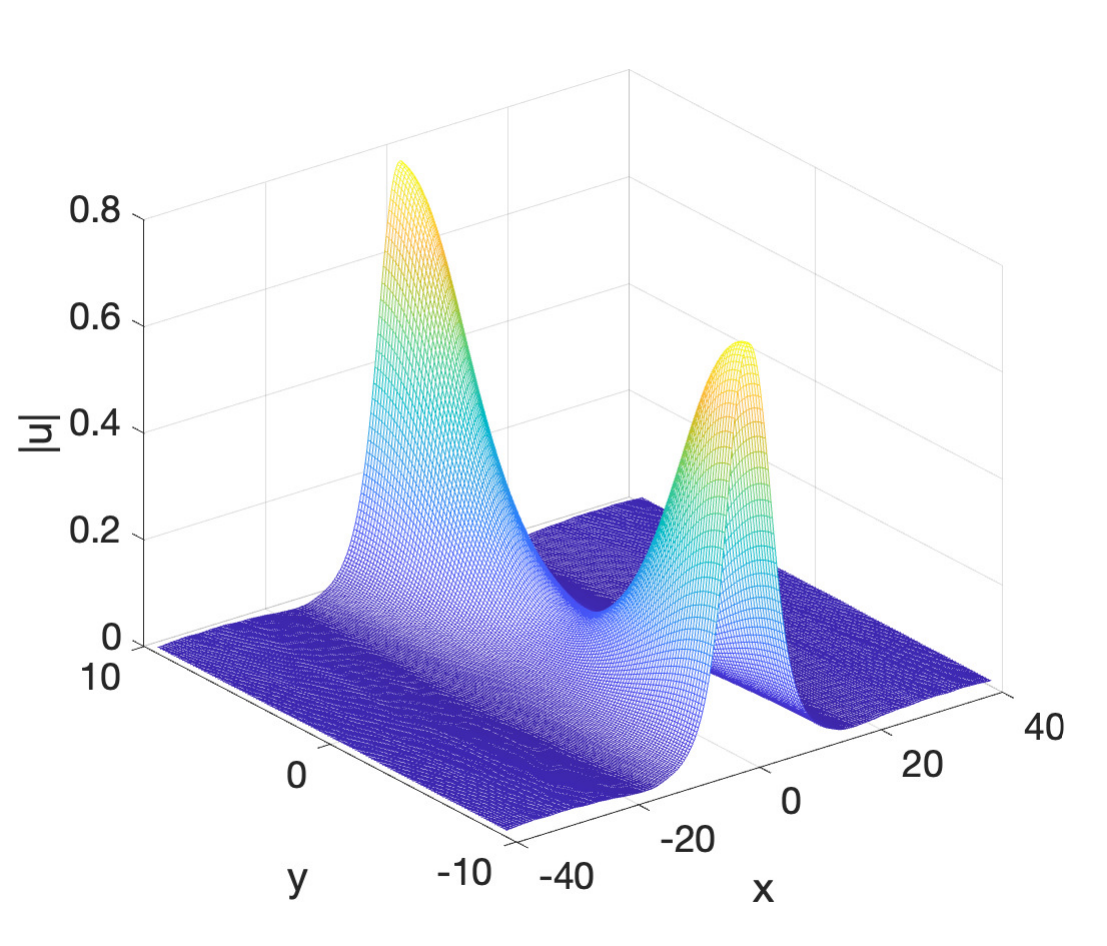}
 \caption{Solution  $|u|$ of (\ref{nls}) at time $t=500$ for initial 
 data (\ref{pert}) with $\omega=0.1$, and 
 $\lambda= - 0.05$ and $L_{y}=3$. }
 \label{figom01_3m}
\end{figure}

If we use the same numerical parameters as before but with $\lambda=-0.05$, we find a 
similar behavior. This time the peak of the final state appears on the 
boundary of the computational domain $\mathbb T_{L_y}$, see Fig.~\ref{figom01_3m}. 
This is consistent with the blow-up at two separate points observed in the solution of the cubic NLS, cf. Fig. \ref{pertcubm}.

\begin{figure}[htb!]
  \includegraphics[width=0.49\textwidth]{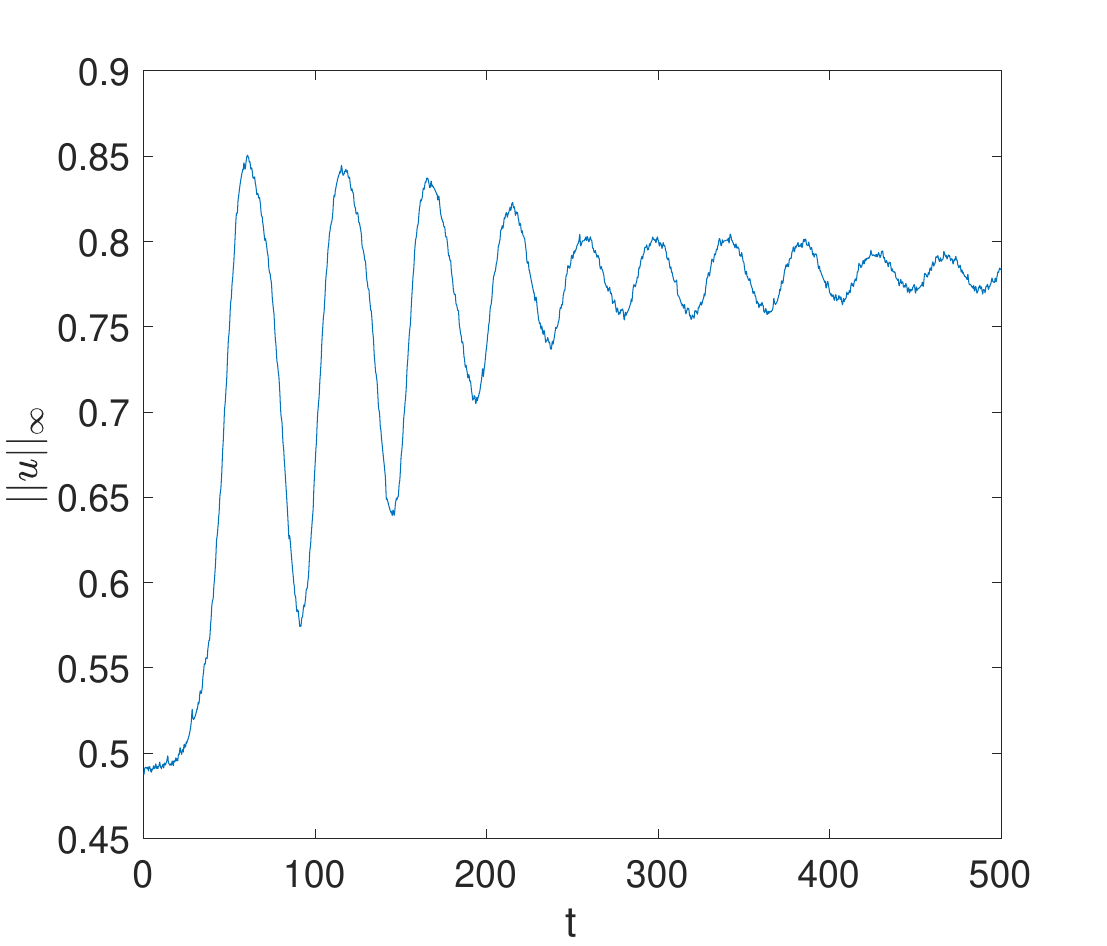}
  \includegraphics[width=0.49\textwidth]{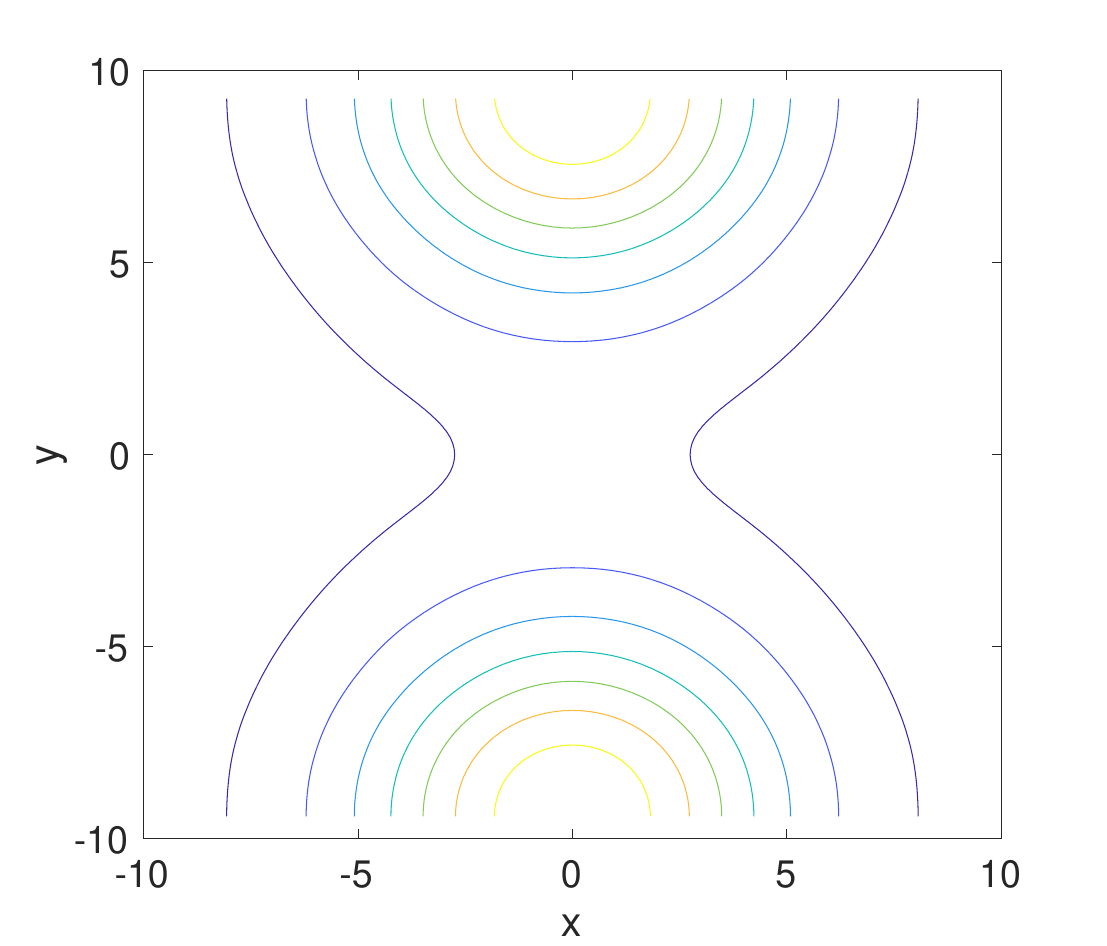}
 \caption{Left: Time-evolution of the $L^{\infty}$ norm of the solution to (\ref{nls}) for the initial 
 data (\ref{pert}) with $\omega=0.1$, $\lambda= - 0.05$ and $L_y =3$. Right: 
 a contour plot of $|u|$ at the final 
 time. }
 \label{figom01_3minf}
\end{figure}

The $L^{\infty}$ norm of the solution $u$ as a function of time is shown on the left of 
Fig.~\ref{figom01_3minf}. Its behavior is qualitatively similar to the previous case.

The contour plot shown on the right of the same figure again indicates radial symmetry of the final state, but this time with a center close to boundary of 
the fundamental cell in $y$. 

%%%%%%%%%%%%%%%%%%%%%%%%

\section{Dependence of the instability on the frequency} 

\begin{figure}[htb!]
  \includegraphics[width=0.49\textwidth]{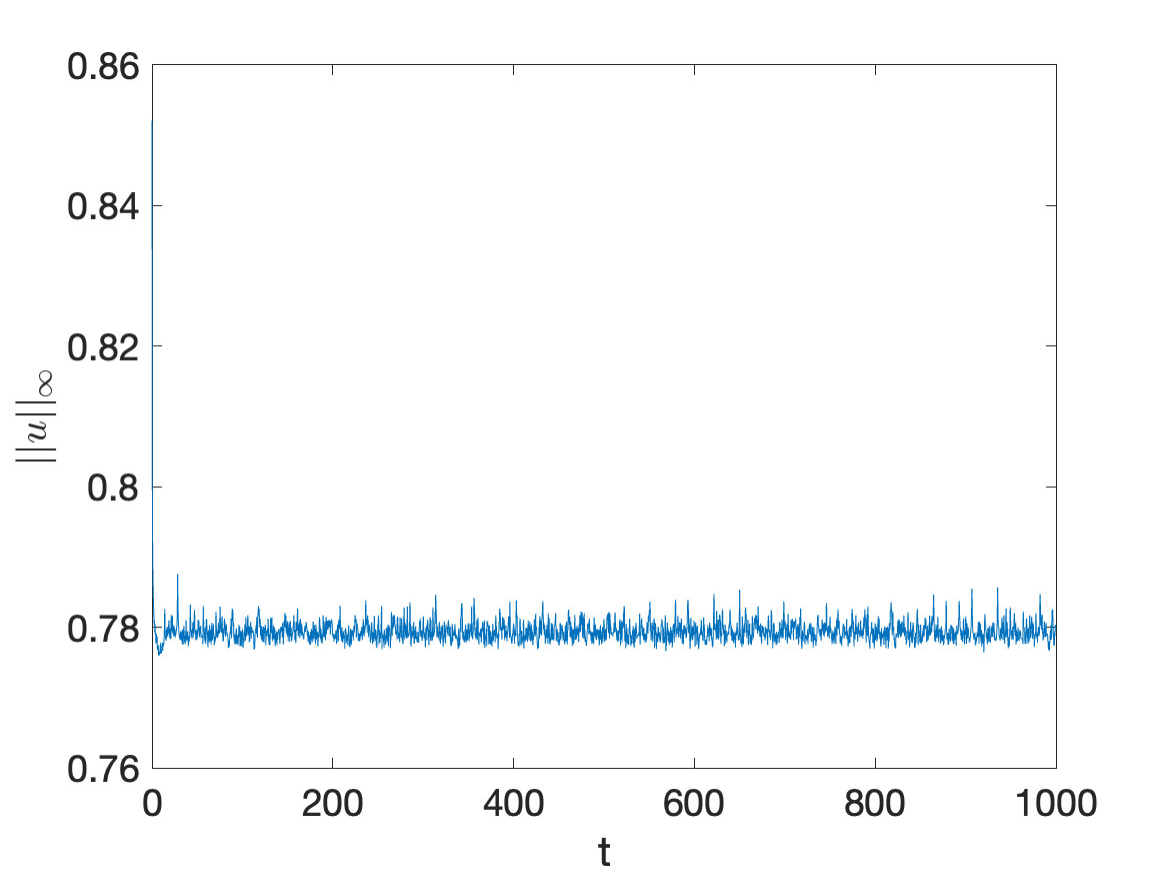}
  \includegraphics[width=0.49\textwidth]{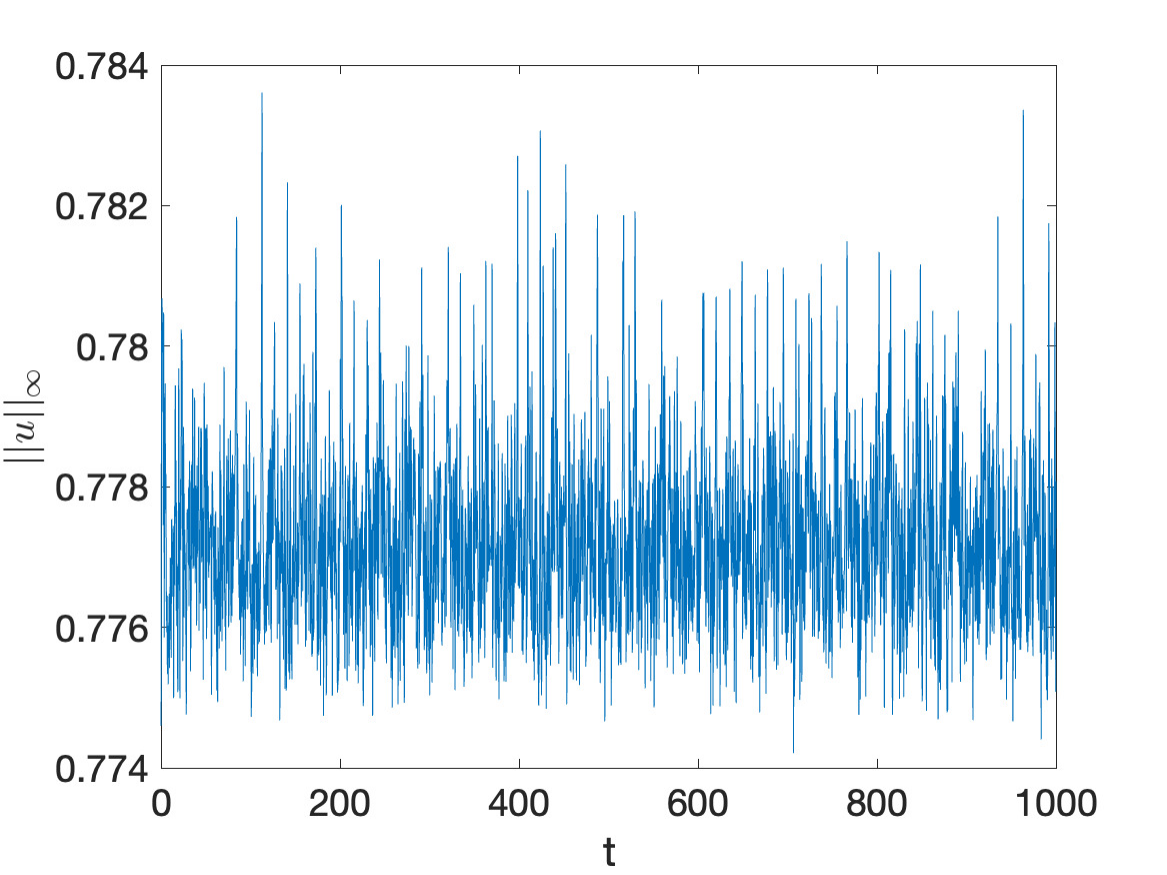}
 \caption{Time-dependence of the $L^{\infty}$ norm of the solution to \eqref{nls} for initial 
 data \eqref{pert} with $\omega=0.18$, $\lambda= \pm 0.077$ and $L_{y}=3$ (on the 
 left for the $+$ sign and on the right for the $-$ sign). }
 \label{figom018_3inf}
\end{figure}

Numerically it is difficult to identify the critical value of $L_{y}$ for which the line solitary wave becomes unstable at a 
given frequency $\omega \in (0, \frac{3}{16})$. Indeed we expect the instability to appear at later times as $\omega$ becomes larger. This can be inferred
 from the heuristic considerations in Section \ref{sec:heur} which imply the need
of a sufficiently large $M(u_0)$ in the considered fundamental cell to allow for
the formation of a lump solitary wave $Q_\omega$. Since the mass in 2D grows 
much faster than the mass in 1D as  $\omega \to \frac{3}{16}$ (cf. Figs. 1 and 3), the line solitary 
wave for increasing values of $\omega$ is expected to remain stable for ever growing values of $L_{y}$. We shall briefly  
illustrate this for the value $\omega = 0.18 \approx \frac{3}{16}$: 

For initial data of the form \eqref{pert} a perturbation of the order of $10\%$ of $\| \phi_{\omega=0.18}\|_{L^\infty}$ yields $\lambda \approx \pm 0.077$. 
Using the same numerical parameters as in Section \ref{sec:unstab}, the perturbed line soliton in this case appears to remain stable even for $L_y=3$, i.e., the torus length for which 
$\phi_{\omega=0.1}$ clearly displayed instability. 
We show $\| u(t, \cdot, \cdot) \|_{L^{\infty}}$ as a function of time for both signs of $\lambda$ in 
Fig.~\ref{figom018_3inf}. In both cases there is no indication of an instability even if one computes until very long times $t\in [0, 1000]$.  
Note that there is considerable amount of radiation (which cannot escape our numerical domain) in the case with the $-$ sign.  

\begin{figure}[htb!]
  \includegraphics[width=0.49\textwidth]{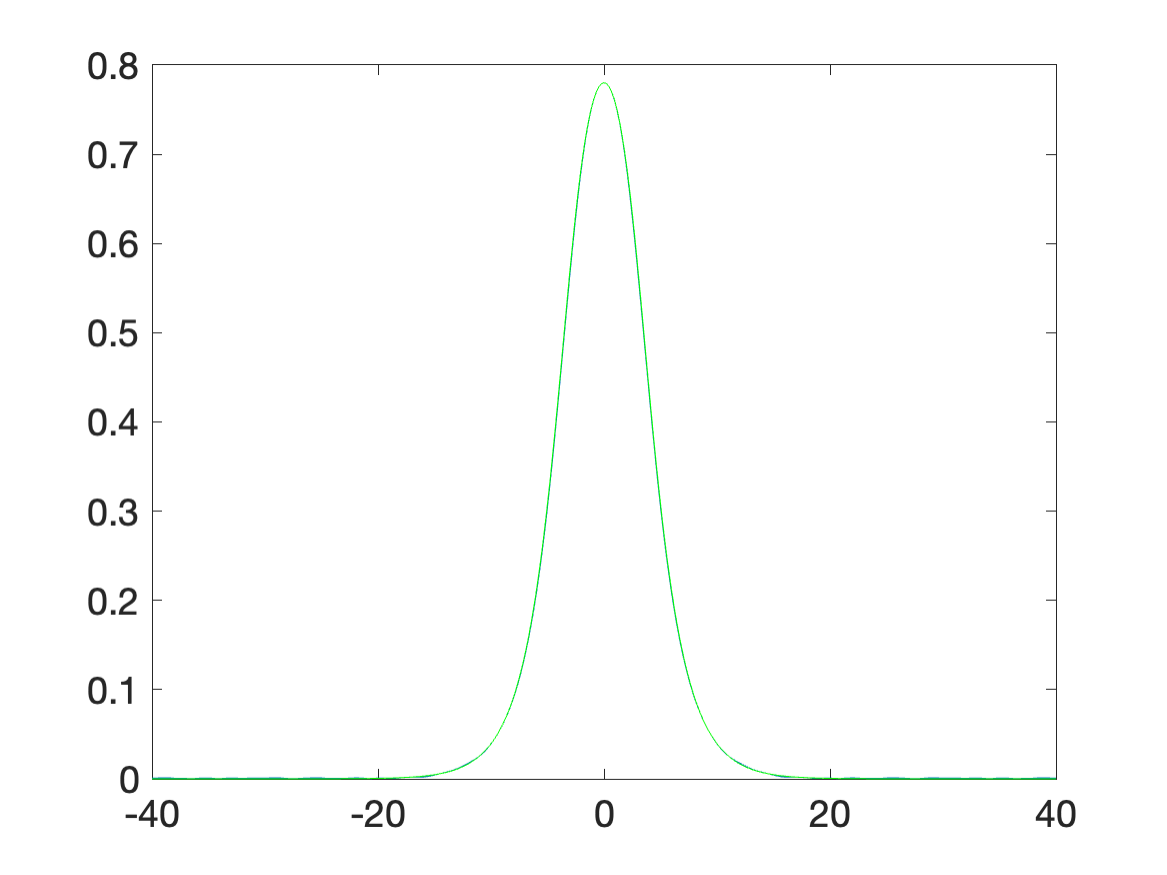}
  \includegraphics[width=0.49\textwidth]{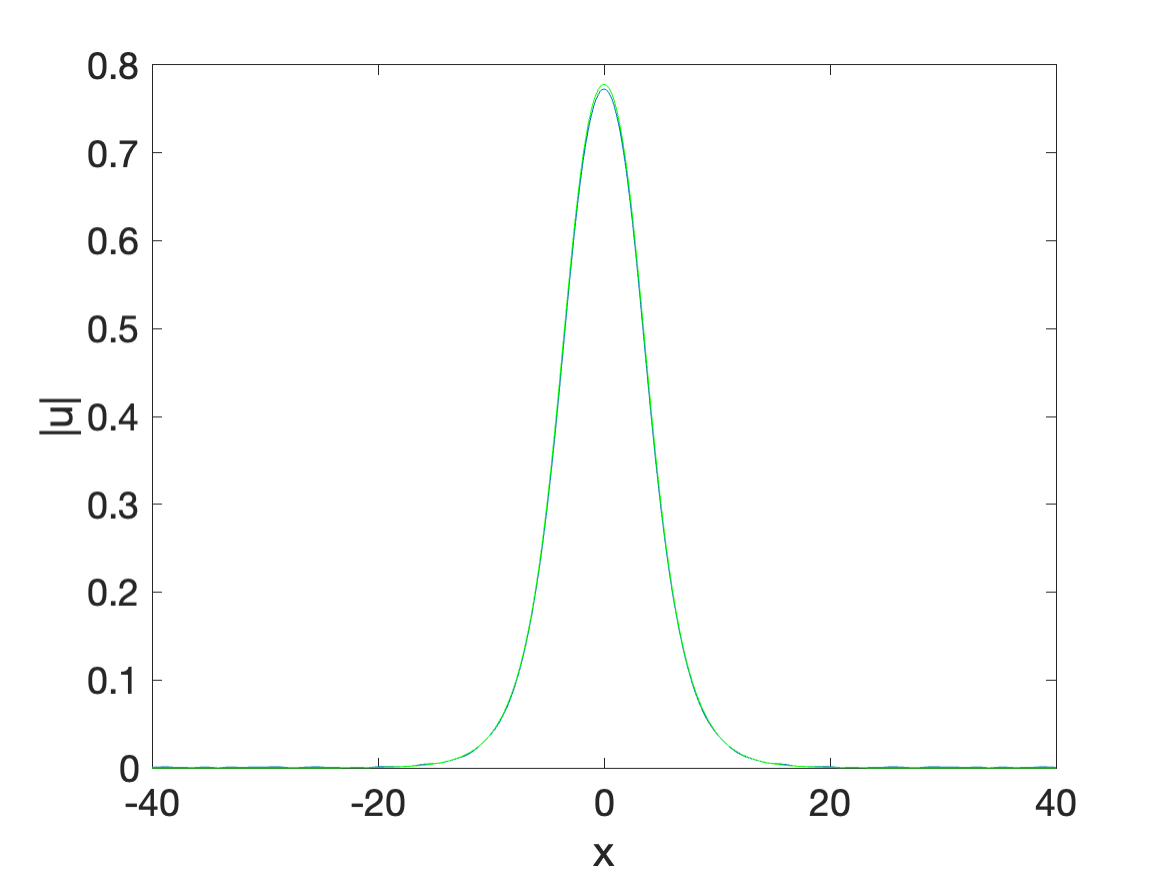}
 \caption{The modulus of the solution to \eqref{nls} for initial date \eqref{pert} with $\omega=0.18$,  
 $\lambda = \pm 0.077$ and $L_{y}=3$ plotted at the final time $t=1000$ along the $y$-axis (in blue),
 together with a fitted solitary wave $\phi_{\omega_\ast}$ in green (on the 
 left for the $+$ sign and on the right for the $-$ sign). }
 \label{figom018_3fit}
\end{figure}

The solution itself is observed to remain close to a line solitary wave $\phi_\omega$ with $\omega=\omega_\ast \not = 0.18$. If we fit the last recorded solution to a line 
solitary wave as described above, we find values of $\omega_\ast =0.1809$ and 
$\omega_\ast =0.1805$ for $\lambda = \pm 0.077$, respectively. 
The solution on the $y$ axis in both cases together with the fitted 
line solitary wave can be seen in Fig.~\ref{figom018_3fit}. 

\begin{figure}[htb!]
  \includegraphics[width=0.49\textwidth]{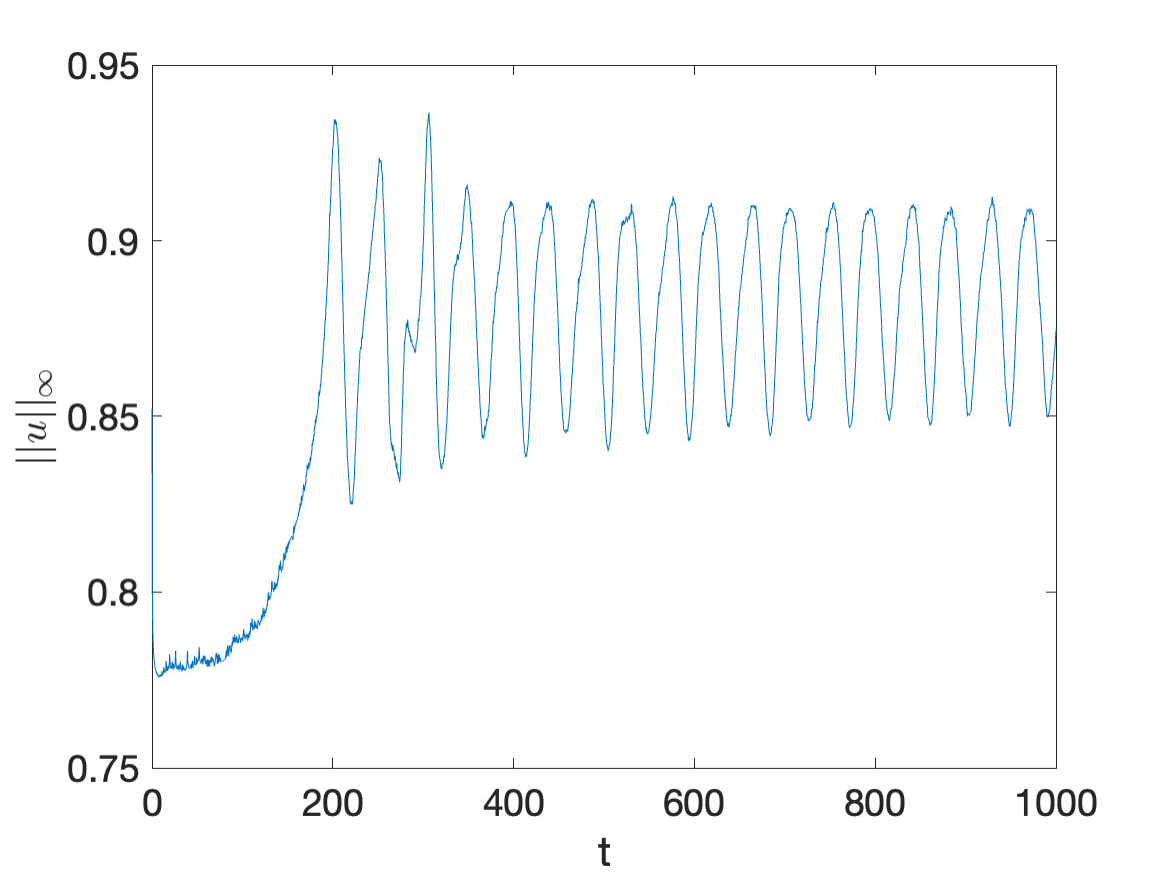}
  \includegraphics[width=0.49\textwidth]{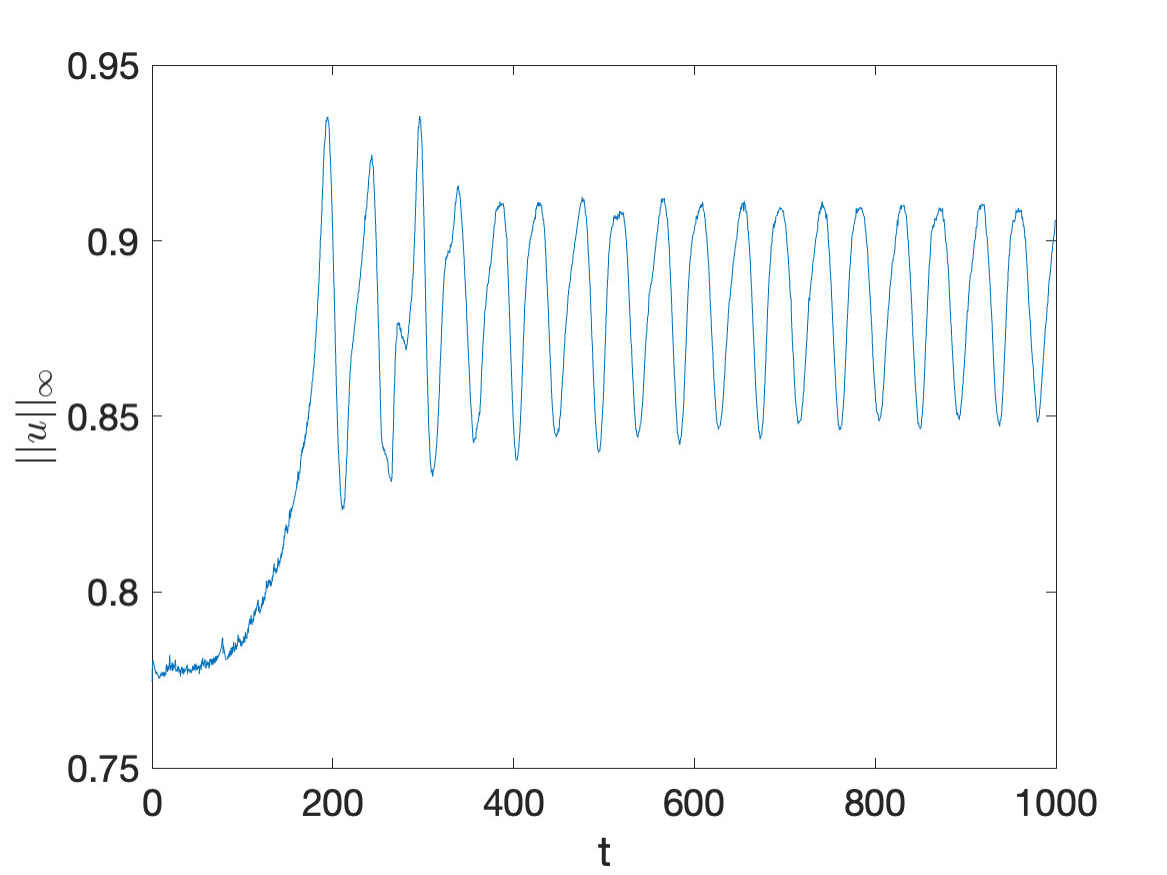}
 \caption{$L^{\infty}$ norm of the solution to \eqref{nls} for  initial 
 data (\ref{pert}) with $\omega=0.18$, $L_{y}=3$ and 
 $\lambda = \pm 0.077$ (on the 
 left for the $+$ sign and on the right for the $-$ sign). }
 \label{figom018_5inf}
\end{figure}

We now consider the same initial data on an even larger torus with length $L_{y}=5$. The numerical parameters in this case are 
$L_{x}=150$, $N_{x}=2^{12}$, 
$N_{y}=2^{7}$, with $N_{t}=2\cdot 10^4$ time steps for $t\in [0,1000]$. 
The $L^{\infty}$ norm of $u$ as a function of time for both signs of $\lambda$ is shown in 
Fig.~\ref{figom018_5inf}. It can be seen that the norm initially grows quite rapidly but then starts to oscillate (with some large period) around a final state presumed to 
be a stable lump solitary wave $Q_{\omega}$. Notice that the amplitude of the oscillations within $\| u (t, \cdot, \cdot) \|_{L^\infty}$ decreases only very slowly as $t\to +\infty$. 
Indeed, even if we compute up to a final time $t=2000$ and take $L_x$ even larger, the picture does not change considerably. 
We expect that a marked decrease of the amplitude of these oscillations (as previously observed in 
Fig.~\ref{figom01_3minf}) can only be observed on time scales much 
larger than those accessible with the direct numerical approach used in the present paper.

\begin{figure}[htb!]
  \includegraphics[width=0.49\textwidth]{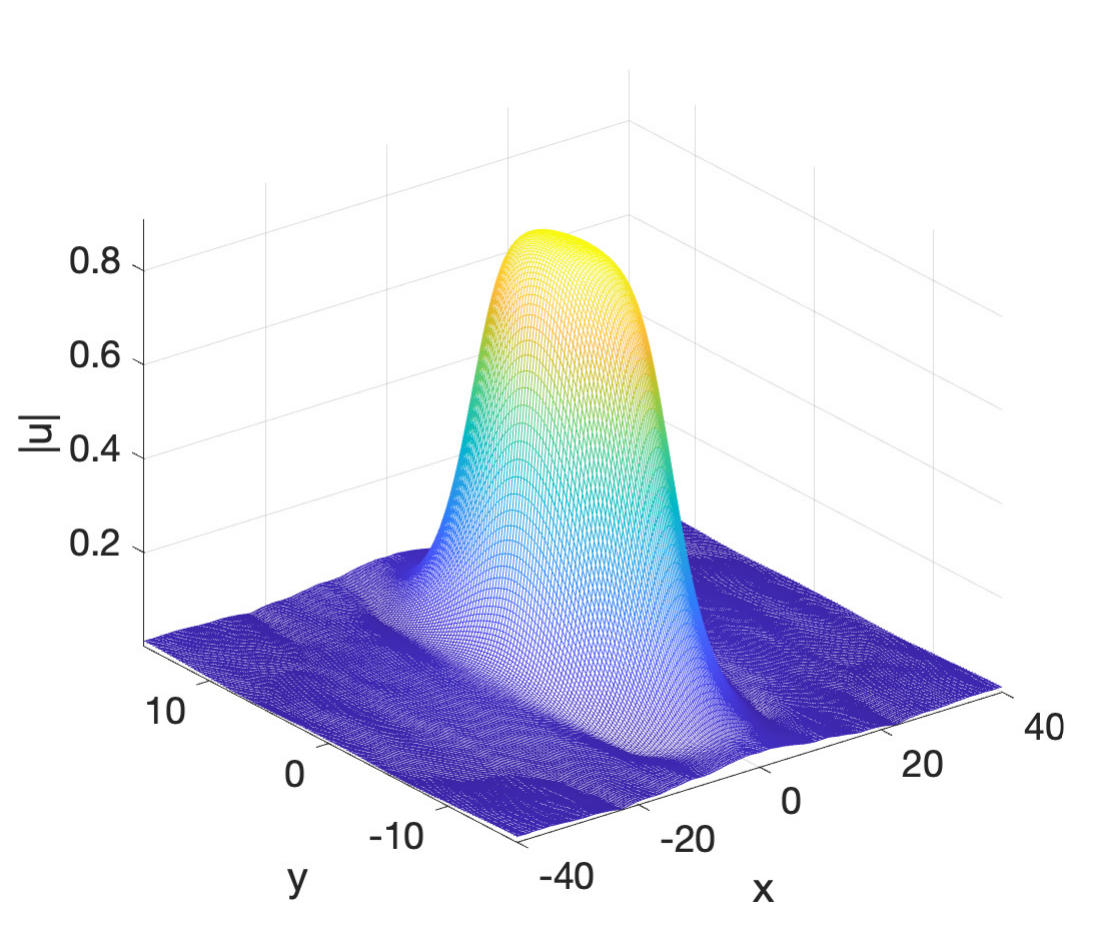}
  \includegraphics[width=0.49\textwidth]{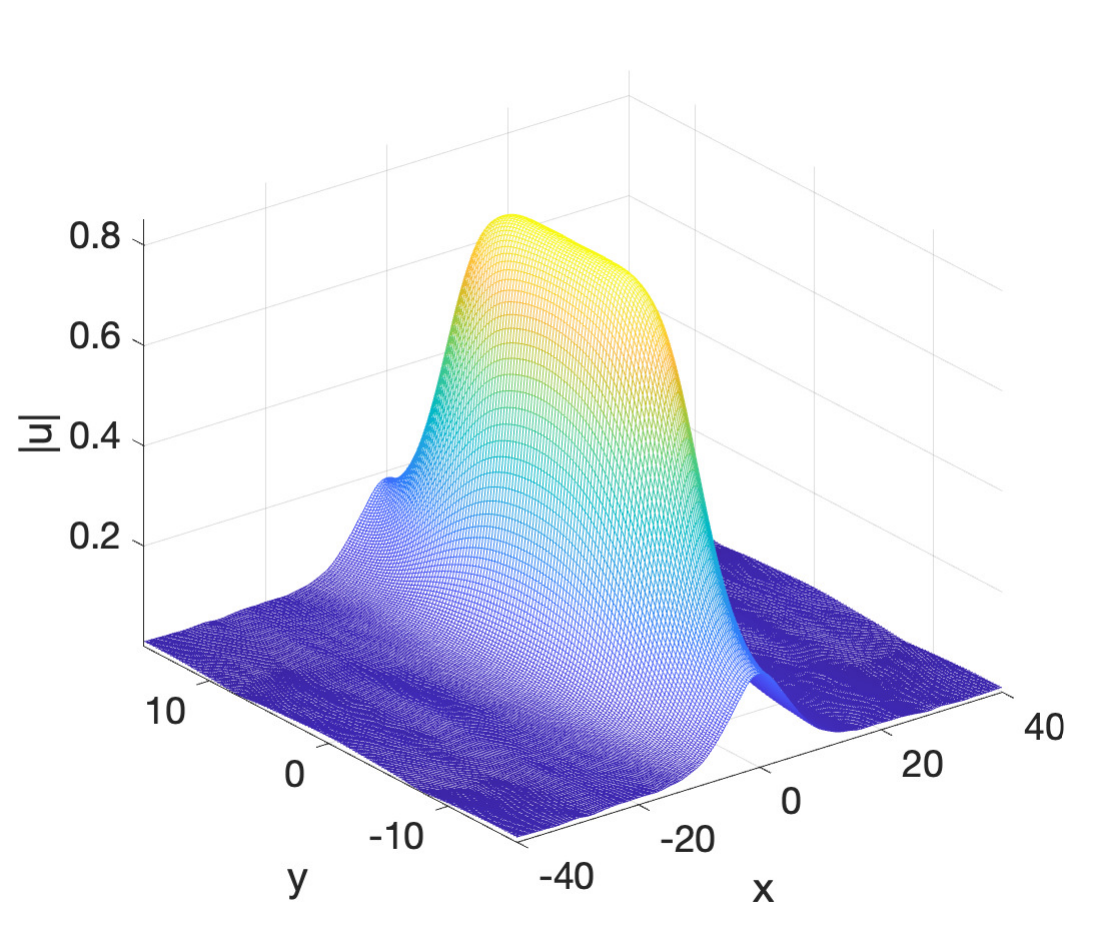}\\
    \includegraphics[width=0.49\textwidth]{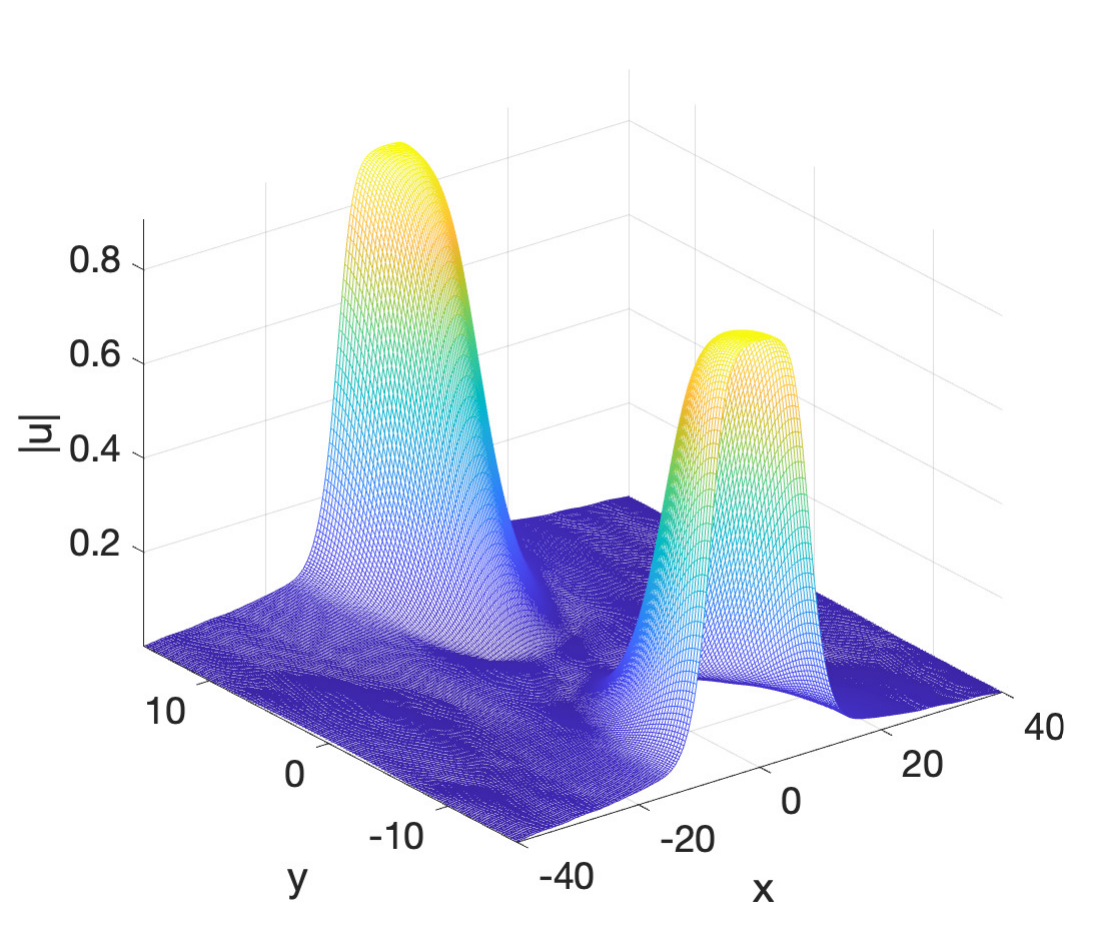}
  \includegraphics[width=0.49\textwidth]{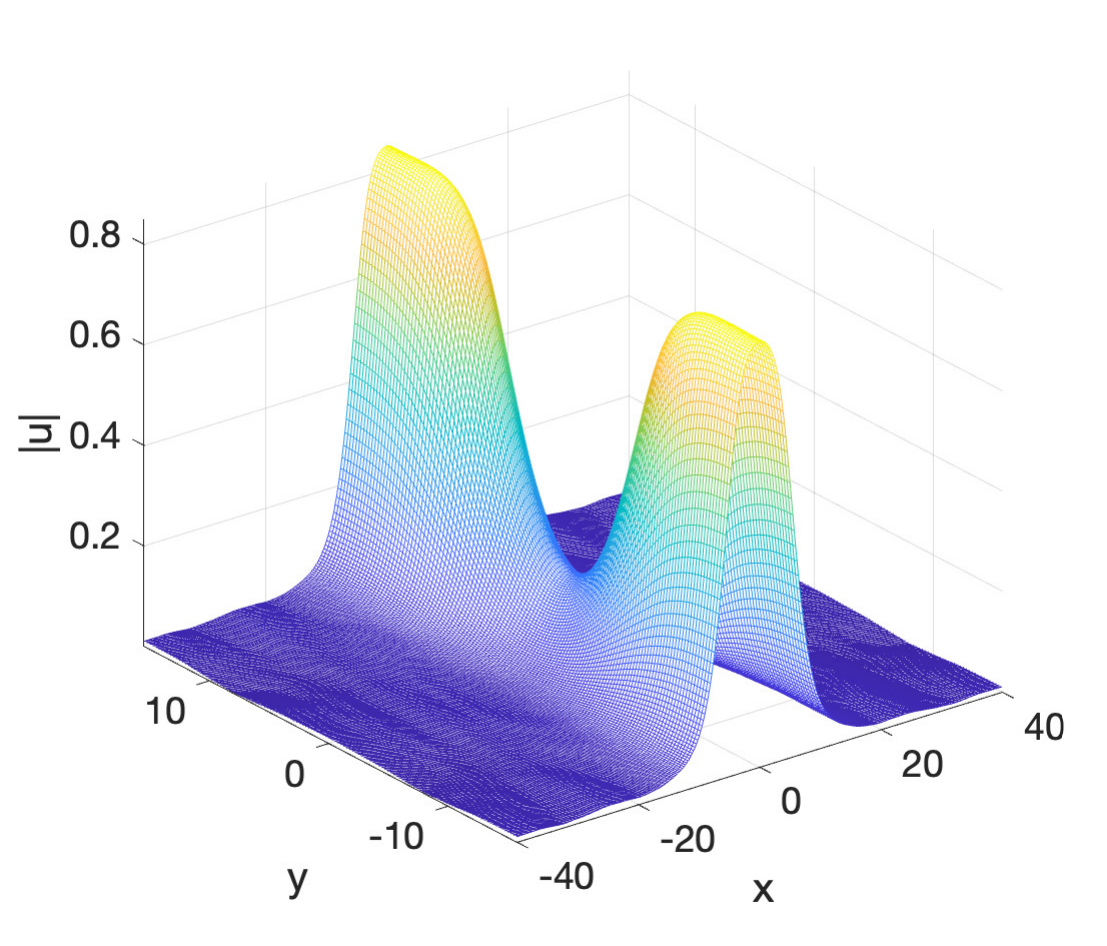}
 \caption{Solution $|u|$ to \eqref{nls} for the initial 
 data \eqref{pert} with $\omega=0.18$,  
 $\lambda = \pm 0.077$ and $L_{y}=3$: in the upper row for the $+$ sign 
 at $t=970$ and $t=990$, in the lower row for the $-$ sign at 
 $t=955$ and $t=980$.}
 \label{figom018_5t1000}
\end{figure}

In Fig.~\ref{figom018_5t1000} we show the states between which the 
solution $u(t, \cdot, \cdot)$ oscillates for $t\in [0,1000]$. In the upper row one can see the 
situation corresponding to the choice $\lambda = 0.077$. The solution visibly oscillates between an 
elongated structure, reminiscent of the initial line solitary wave, and a more peaked structure which appears to include a part of the lump solitary wave. 
The situation is similar for the case $\lambda = -0.077$, shown in the 
lower row of the same figure, except that the peaked structures appear at the boundary of the numerical domain.

%%%%%%%%%%%%%%%%%%%%%%%%%%%%%%%%%%%

\section{Periodic deformations}\label{sec:perper}

In this final section we shall study a class of perturbations which is only admissible on wave guide domains $\mathbb R_x \times \mathbb T_y$, namely 
periodic deformations of the line solitary wave $\phi_\omega$. To be more precise, we shall in the following consider initial data given by
\begin{equation}\label{definitial}
	u_{0}(x,y) = \phi_{\omega}\big(x-\lambda \cos(y)\big), \quad \lambda \ge 0.
\end{equation}
In our numerical simulations below we will always choose $\lambda = 0.8$, which corresponds to a rather large periodic deformation.
To gain more insight on the behavior of the solution, we shall first consider the case of the cubic NLS and then study the corresponding 
situation for the cubic-quintic NLS. 

\subsection{Periodic deformations in the cubic NLS}
We consider the cubic NLS \eqref{cubic} with initial data \eqref{definitial} at $\omega = 0.04$ and $\omega =1$, respectively. 
We use numerical parameters $L_{x}=100$ for $\omega=0.04$ and $L_{x}=10$ for $\omega=1$, together with $L_{y}=2$, $N_{x}=2^{12}$, $N_{y}=2^{7}$ and 
$N_{t}=10^{4}$ time steps for $t\in [0,100]$. 

\begin{figure}[htb!]
  \includegraphics[width=0.49\textwidth]{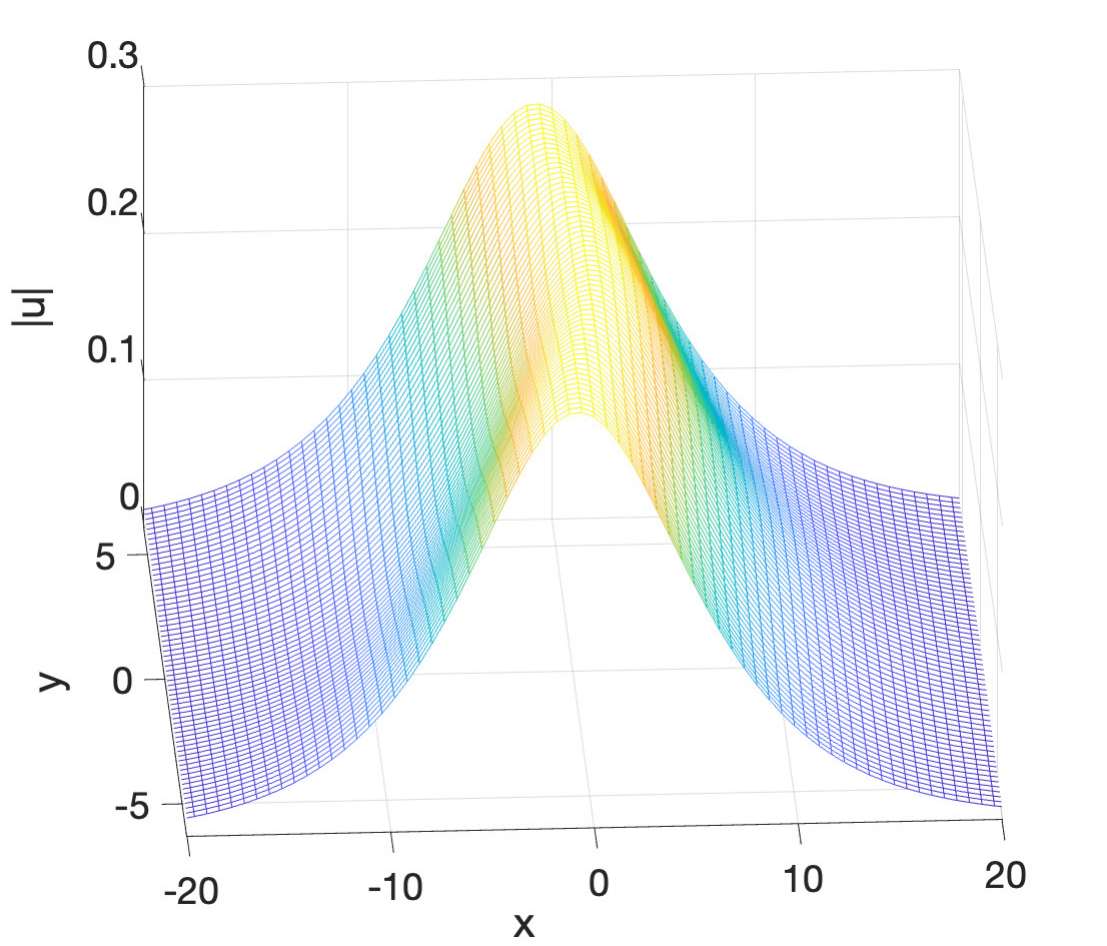}
  \includegraphics[width=0.49\textwidth]{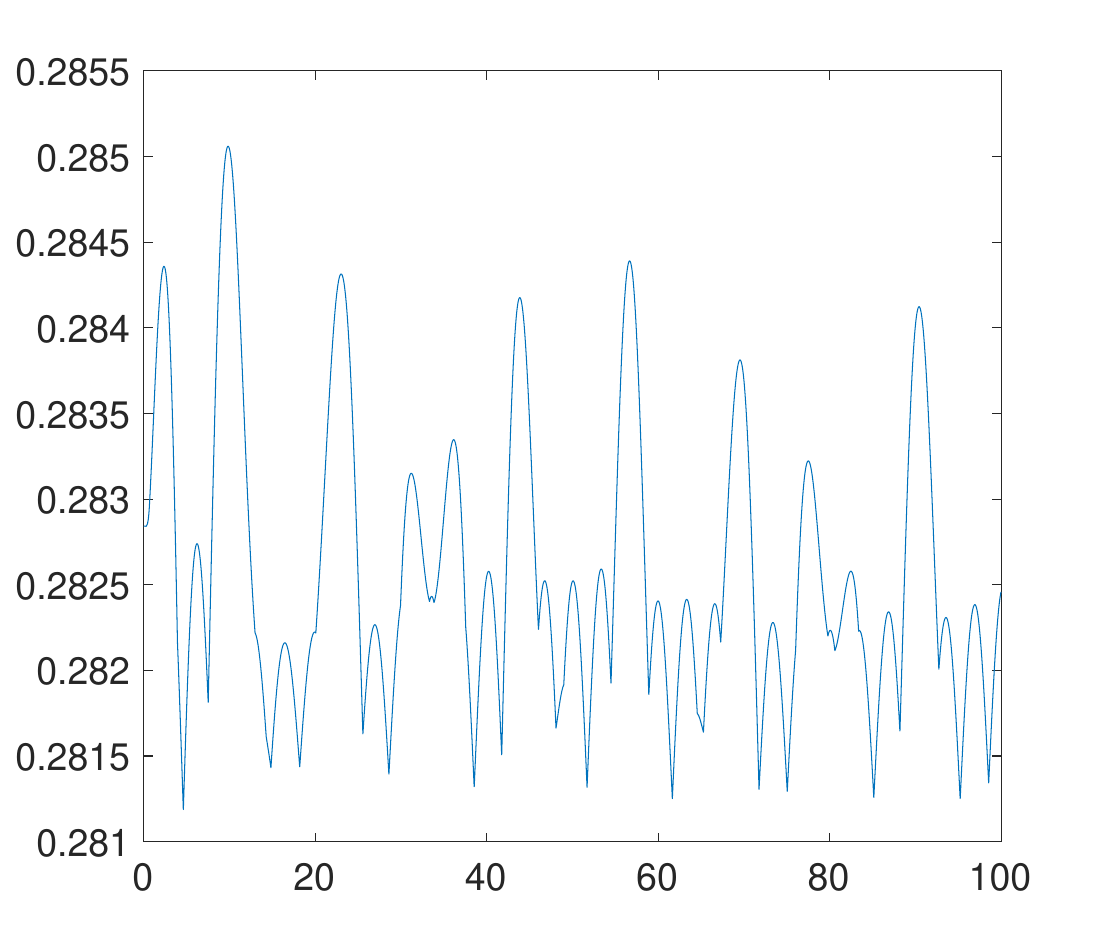}\\
    \includegraphics[width=0.49\textwidth]{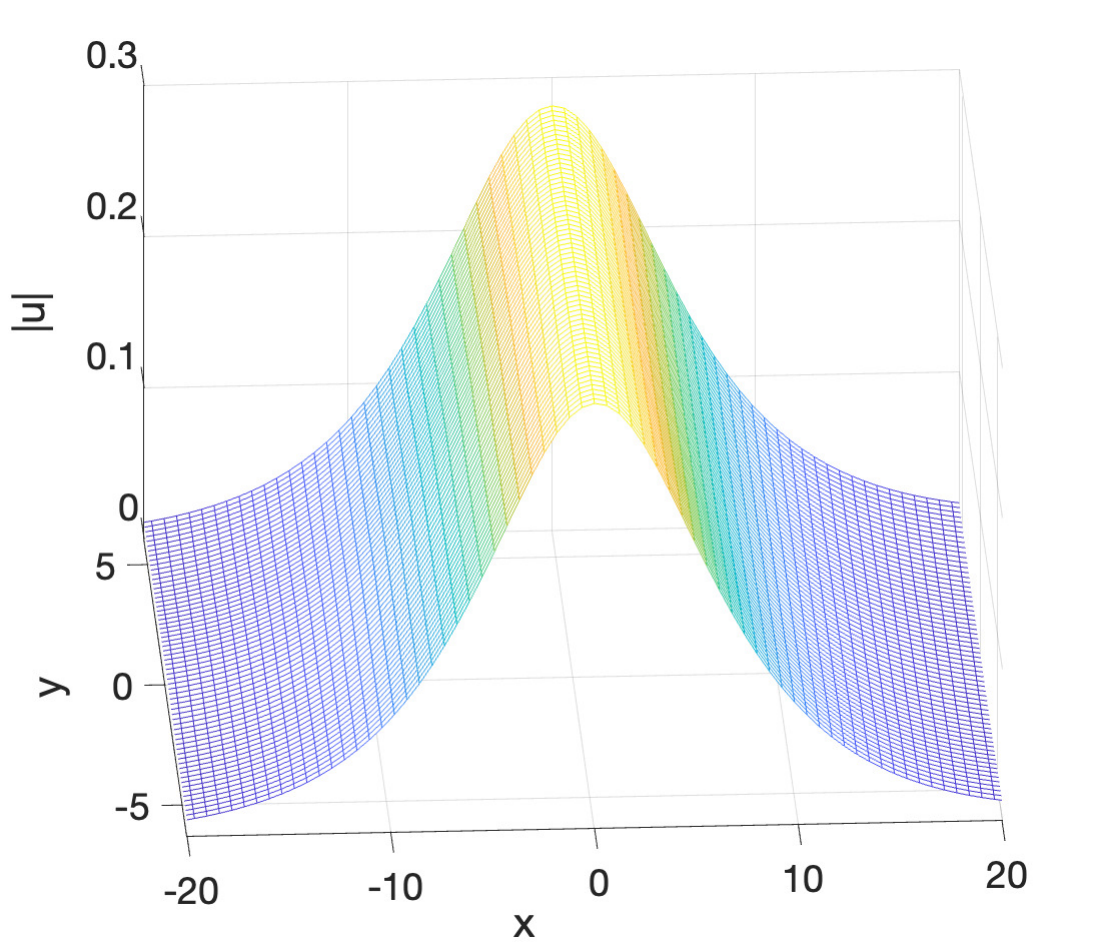}
  \includegraphics[width=0.49\textwidth]{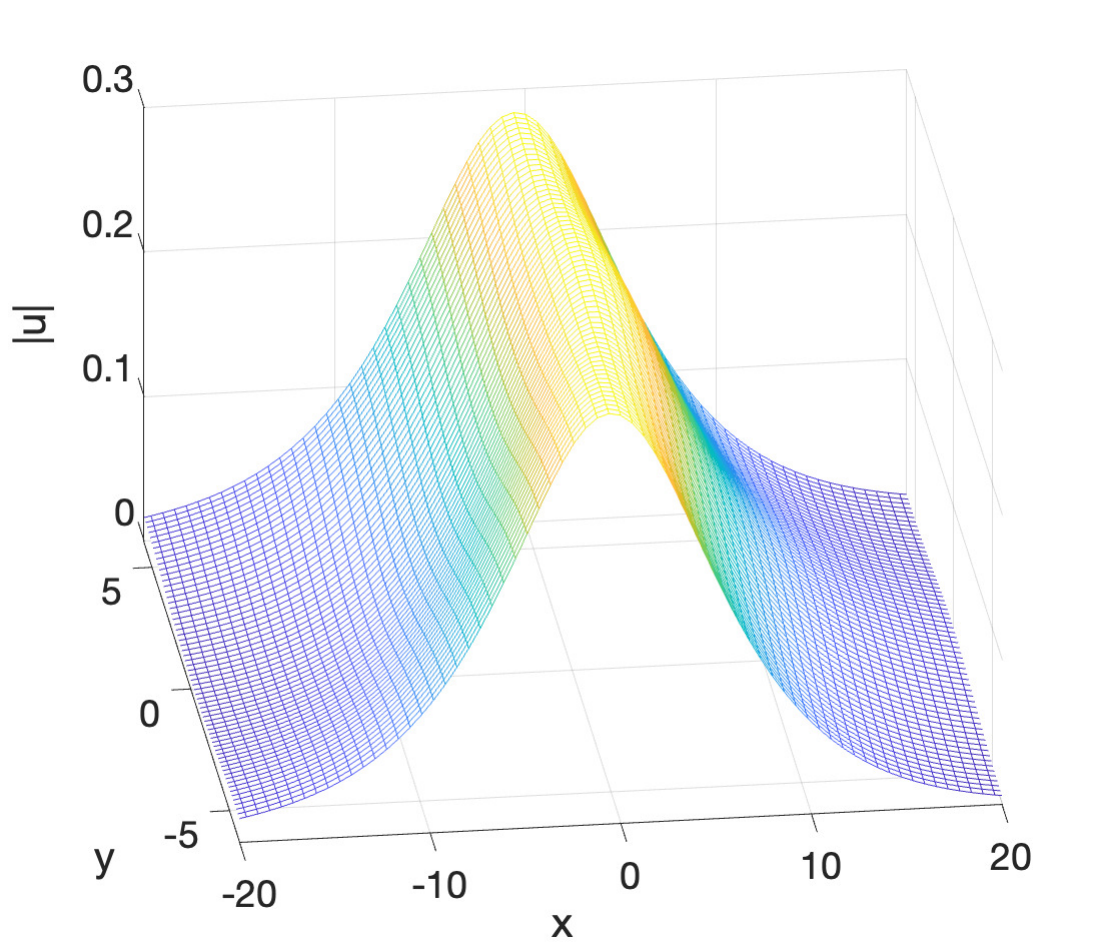}
  \caption{Solution $|u|$ to \eqref{cubic} for initial 
 data \eqref{definitial} with $\omega=0.04$: in the upper row the 
 initial data and the $L^{\infty}$ norm of $u$ as a function of time;  
 in the lower row the solution $|u|$ at times
 $t=85$ and $t=90$.}
 \label{figcubicdefom004}
\end{figure}

First we consider the case $\omega=0.04$, for which we recall that the initial mass $M(\phi_{\omega =0.04}) < M(Q^{\rm cub}_{\omega=0.04})$. 
The corresponding initial data is shown on the left of Fig.~\ref{figcubicdefom004}. The $L^{\infty}$ norm of the solution $u$ as a function of time can be seen on the right of the same figure. 
One observes that there are only some small oscillations within $\| u(t, \cdot , \cdot)\|_{L^\infty}$ and this behavior does not change as $t\to +\infty$. Indeed, it can be seen that the 
solution $u$ oscillates between different 
deformations of the line solitary wave which is the expected behavior in the stable case (see the lower row of the 
same figure). 

\begin{figure}[htb!]
  \includegraphics[width=0.49\textwidth]{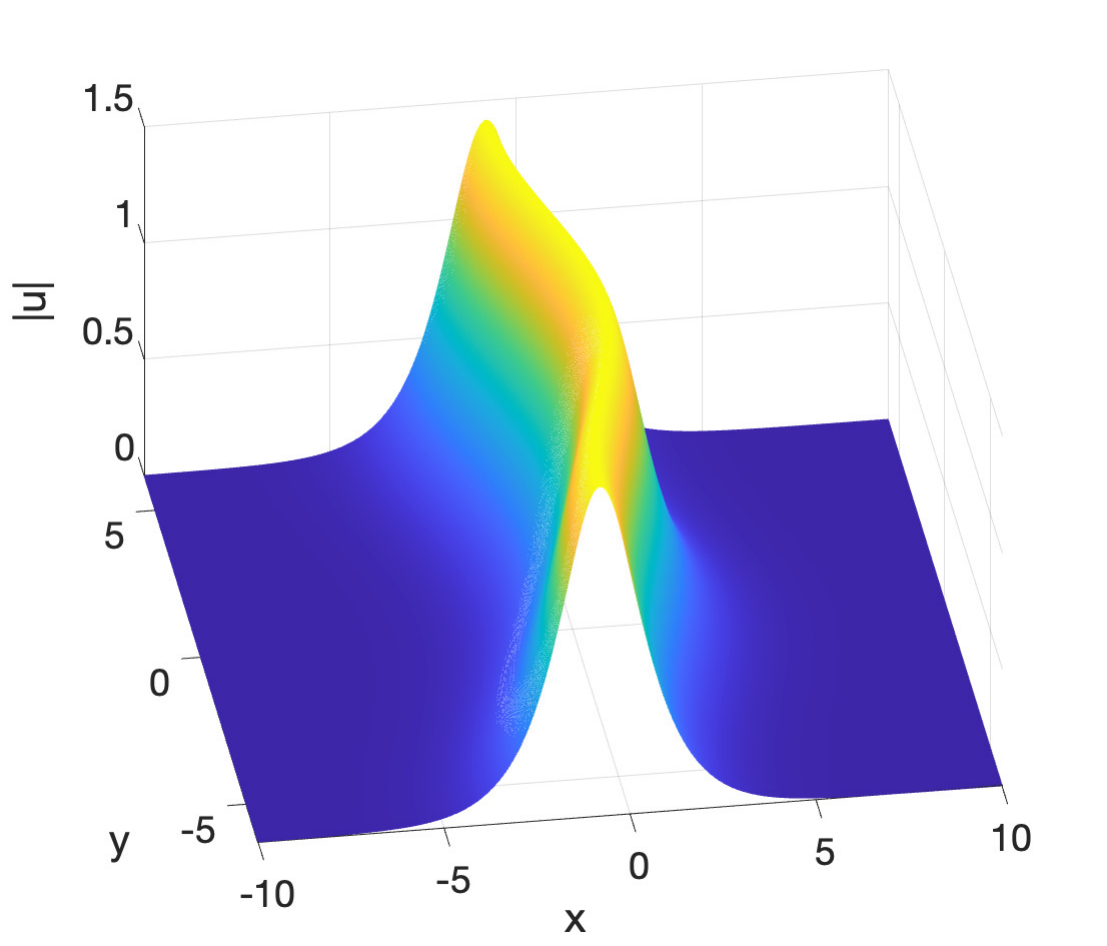}
  \includegraphics[width=0.49\textwidth]{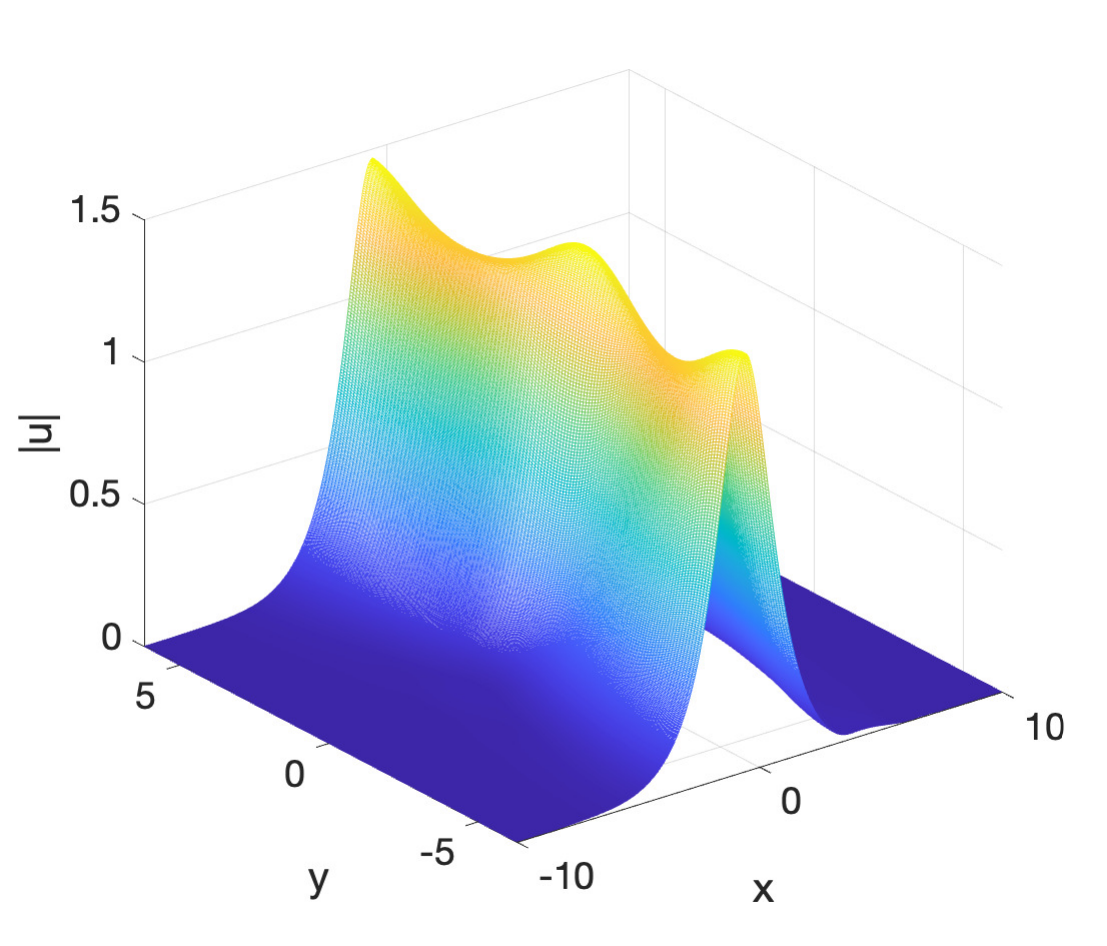}\\
    \includegraphics[width=0.49\textwidth]{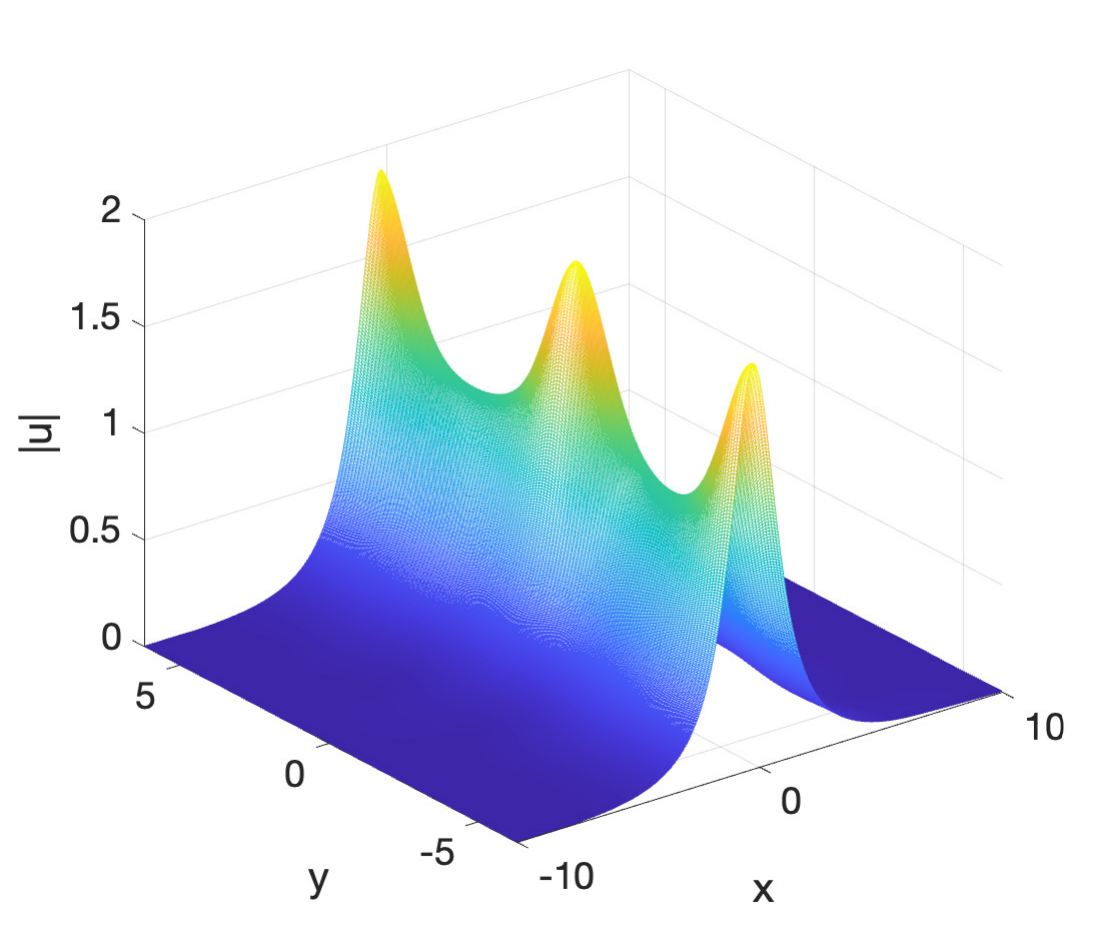}
  \includegraphics[width=0.49\textwidth]{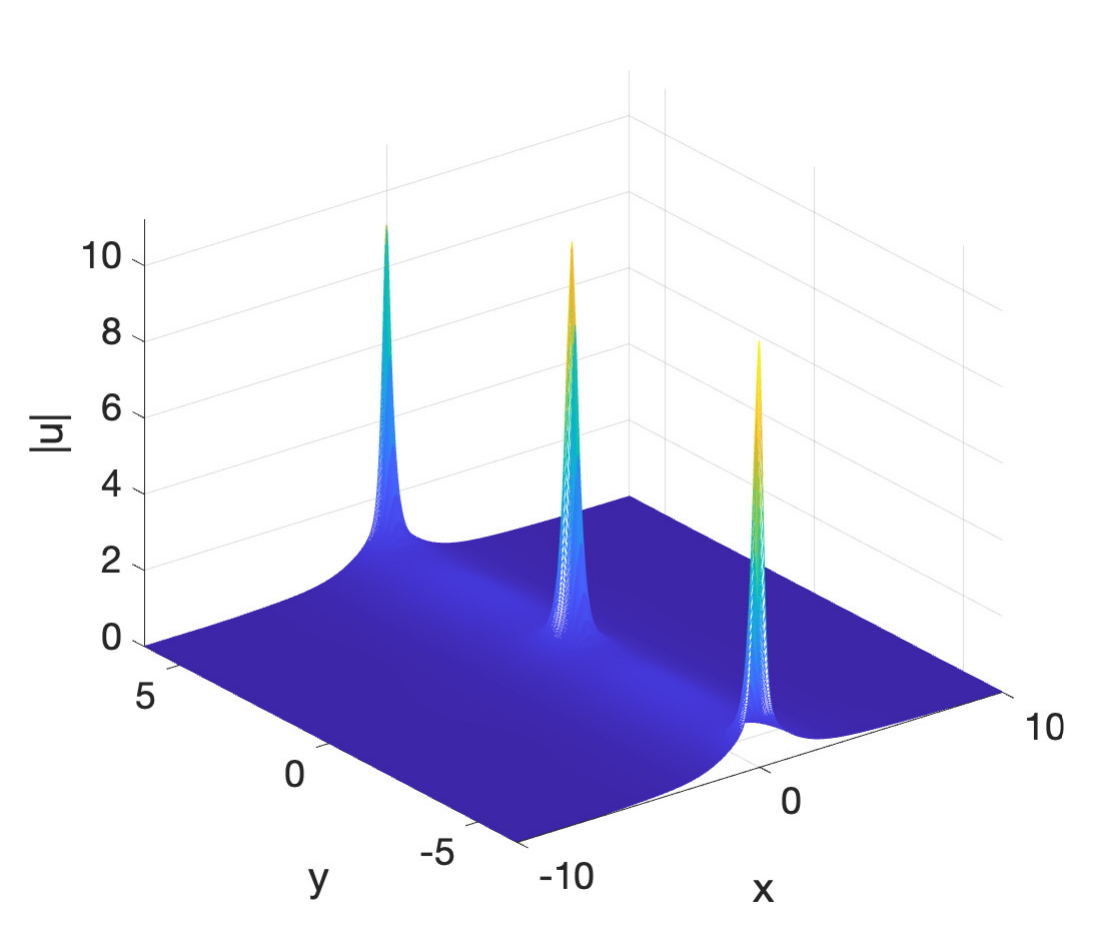}
  \caption{Solution $|u|$ to \eqref{cubic} for the initial 
 data \eqref{definitial} with $\omega=1$ at times $t=0$, $t=1$, $t=2$ 
 and $t=2.75$, respectively.}
 \label{figcubicdefom1}
\end{figure}

In the case $\omega=1$, we recall that $M(\phi_{\omega =1}) > M(Q^{\rm cub}_{\omega=1})$ and thus we expect strong instability due to finite time blow-up. 
The latter is seen to occur simultaneously at three different points (induced by the periodic perturbation) at time $t\approx 2.75$, cf. Fig.~\ref{figcubicdefom1}.

\subsection{Periodic deformations in the cubic-quintic NLS}
For periodically perturbed initial data of the cubic-quintic NLS \eqref{nls} we shall only study the case where $\omega=0.1$, but using different values of the torus length $L_y$. 
We thereby work with numerical parameters $L_{x}=150$, $N_{x}=2^{12}$, $N_{y}=2^{6}$, and $N_{t}=10^{4}$ for $t\in [0, 1000]$.

\begin{figure}[htb!]
  \includegraphics[width=0.49\textwidth]{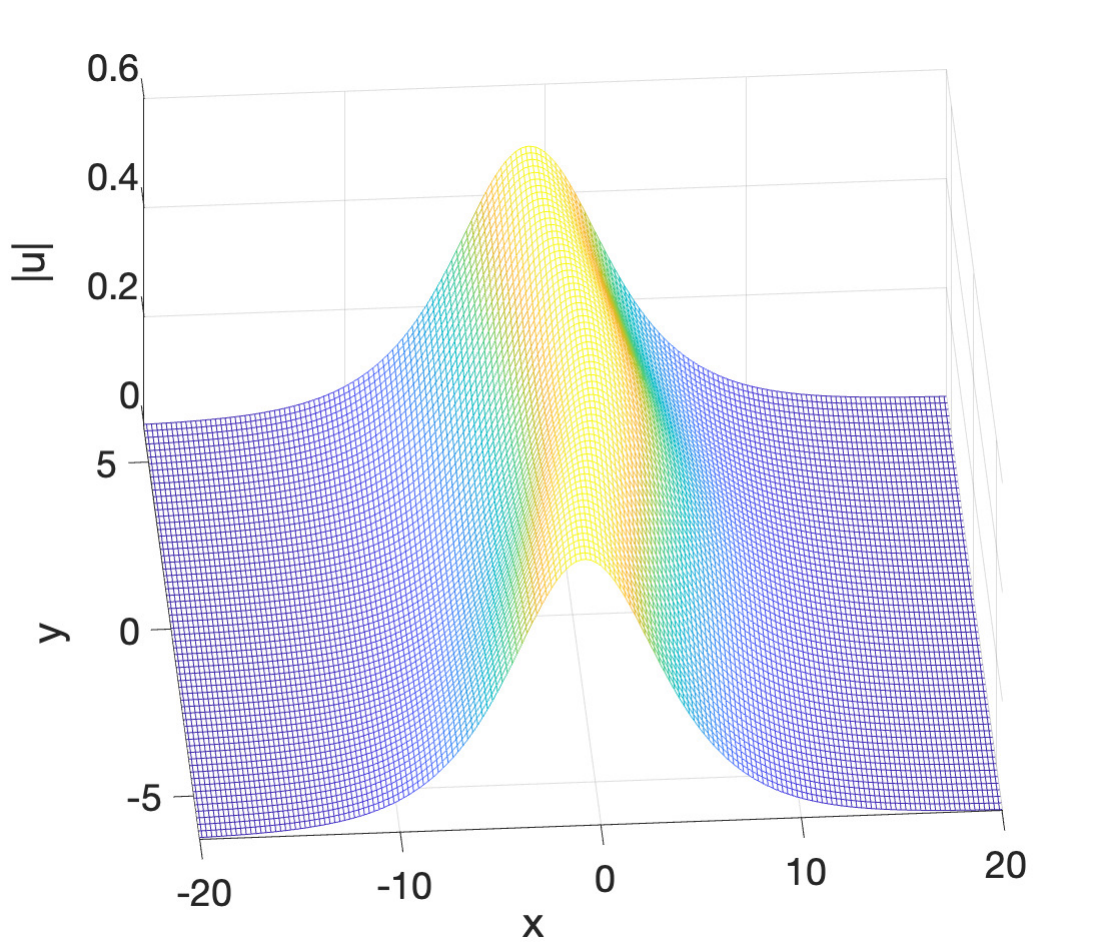}
  \includegraphics[width=0.49\textwidth]{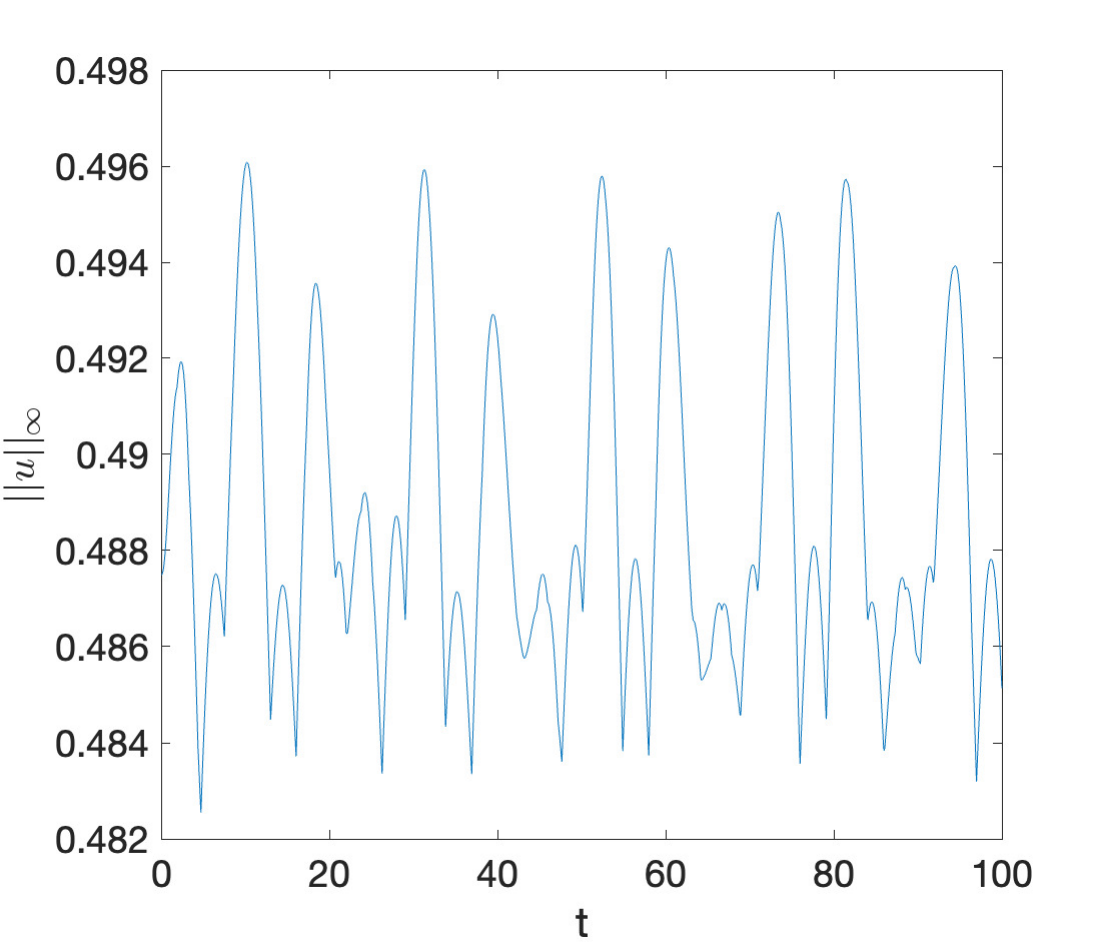}\\
    \includegraphics[width=0.49\textwidth]{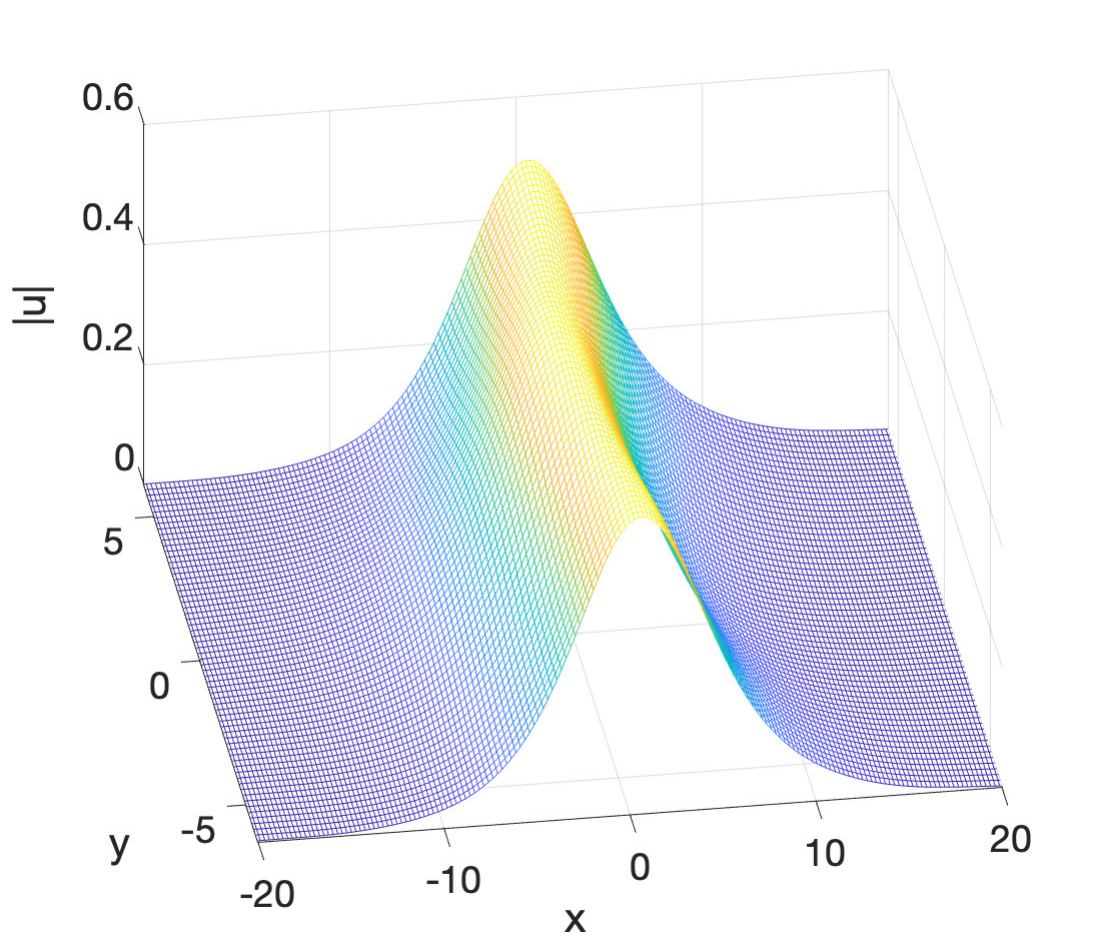}
  \includegraphics[width=0.49\textwidth]{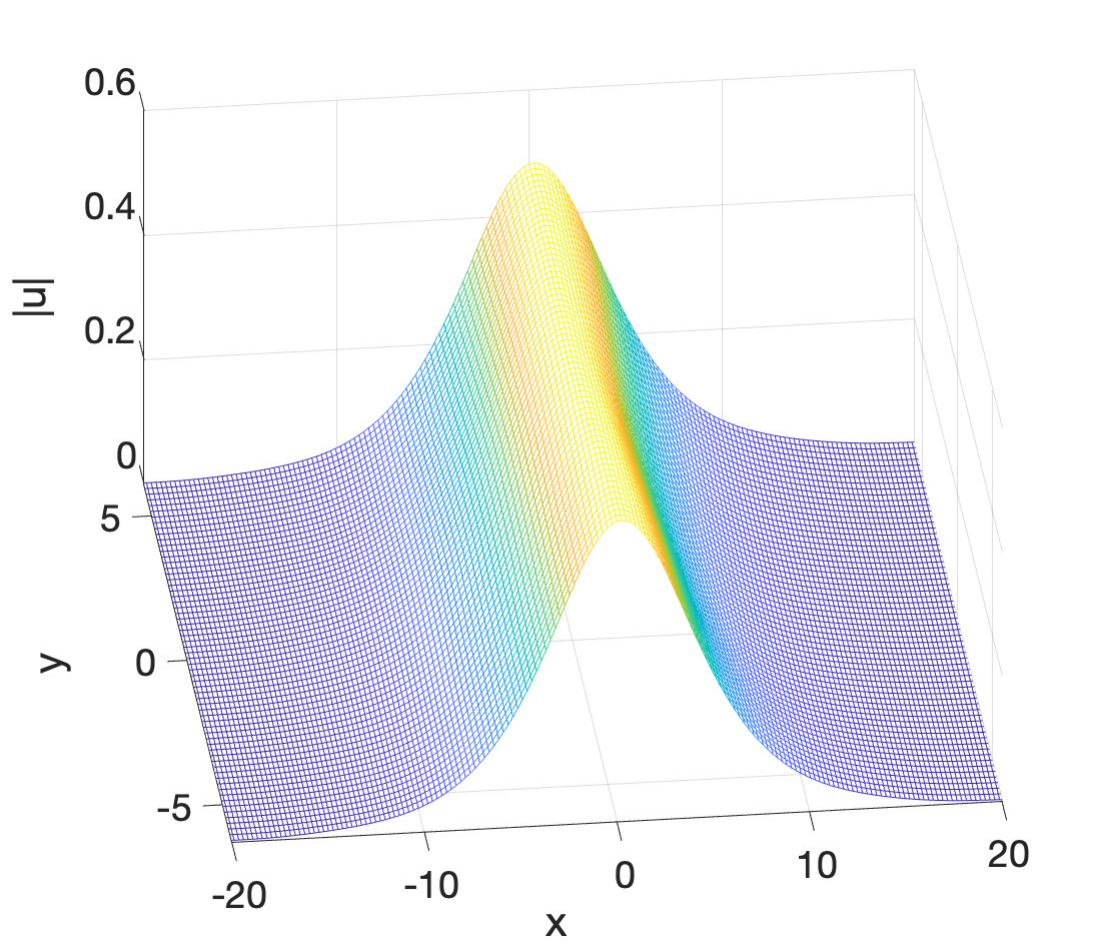}
  \caption{Solution $|u|$ to \eqref{nls} for initial 
 data \eqref{definitial} with $\omega=0.1$ and $L_{y}=2$: in the upper row the  
 the initial data and the $L^{\infty}$ norm of $u$ as a function of time;
 in the lower row the solution $|u|$ for times
 $t=94$ and $t=97$.}
 \label{figcubicquindefLy2}
\end{figure}

We first study the situation with $L_{y}=2$ for which we again expect stability of the periodically perturbed line soliton. That this is indeed the case can be seen in 
Fig.~\ref{figcubicquindefLy2}. Indeed, the behavior of the solution is seen to be qualitatively very similar to the one for the cubic NLS. Compared with Fig. \ref{figcubicdefom004} however, we 
see that the oscillations within $\| u(t, \cdot , \cdot)\|_{L^\infty}$ are larger than in the cubic case. 

\begin{figure}[htb!]
  \includegraphics[width=0.49\textwidth]{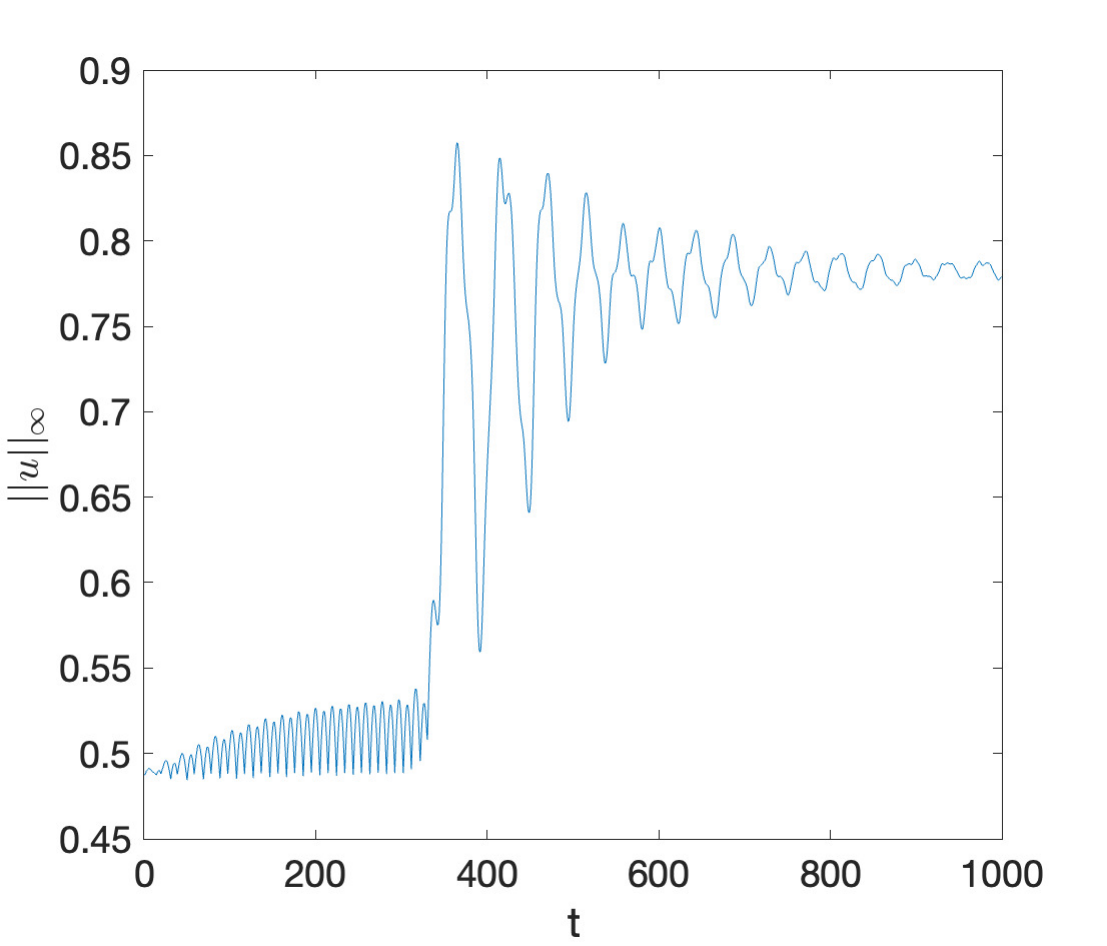}
  \includegraphics[width=0.49\textwidth]{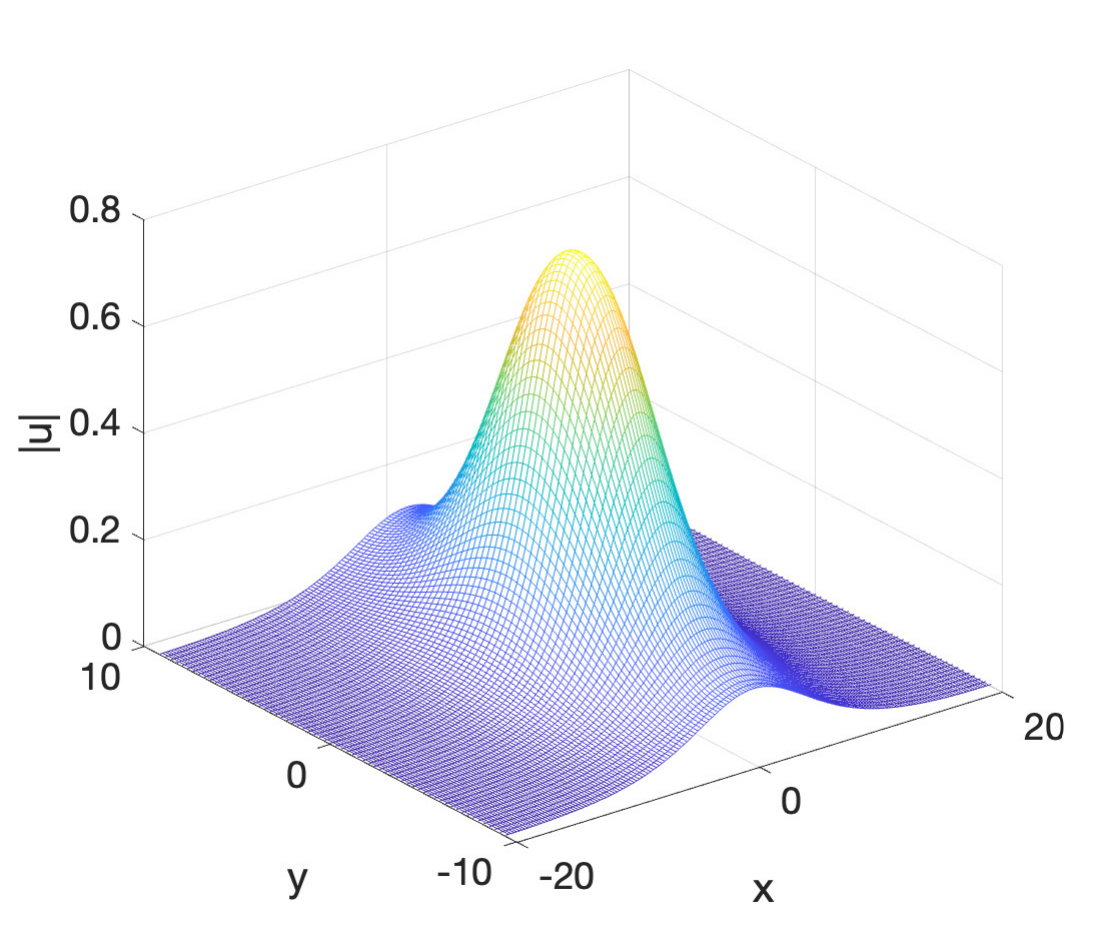}
  \caption{Solution $|u|$ to \eqref{nls} for initial 
 data \eqref{definitial} with $\omega=0.1$ and $L_{y}=3$: on the left 
 the $L^{\infty}$ norm as a function of time and on the right the solution $|u|$ at
 the final time $t=1000$.}
 \label{figcubicquindefLy3}
\end{figure}

Finally, turn to the case $L_{y}=3$ for which we expect instability.  Indeed, it can be see in Fig.~\ref{figcubicquindefLy3} that $\| u(t, \cdot , \cdot)\|_{L^\infty}$, after some 
initial oscillations around the line solitary wave, exhibits a 
dramatic jump at $t\approx 330$. Afterwards one observes large oscillations in time with a decreasing amplitude as $t\to + \infty$. The 
final state of the solution $u$ in this case appears to be again a lump solitary wave $Q_\omega$. 

\begin{remark}
The instability of the periodically perturbed line solitary wave $\phi_\omega$ against the formation of lump solitons $Q_\omega$ is analogous to the situation observed within the KP-I equation, cf. \cite{KMS}. 
The main difference to our case is the fact that KP-I is a unidirectional wave equation and thus its solution cannot oscillate (in time) between the two different types of solitary waves. 
\end{remark}

%%%%%%%%%%%%%%%%%%%%%%%%%%%%%%%%%%%
%%%%%%%%%%%%%%%%%%%%%%%%%%%%%%%%%%%

\bibliographystyle{amsplain}

\end{document}